\documentstyle{amltd}
\begin{document}
\annalsline{157}{2003}
\received{March 30, 1999}
\startingpage{1}
\def\bye{\end{document}}
 \font\tenrm=cmr10
\input amssym.def
\input amssym.tex

\input boxedeps.tex 
\SetepsfEPSFSpecial 
\HideDisplacementBoxes
\def\figin#1#2{
$$
 {\BoxedEPSF{#1.eps scaled
#2}%
}%
$$
\noindent}

\catcode`\@=11
\font\twelvemsb=msbm10 scaled 1100
\font\tenmsb=msbm10
\font\ninemsb=msbm10 scaled 800
\newfam\msbfam
\textfont\msbfam=\twelvemsb  \scriptfont\msbfam=\ninemsb
  \scriptscriptfont\msbfam=\ninemsb
\def\msb@{\hexnumber@\msbfam}
\def\Bbb{\relax\ifmmode\let\next\Bbb@\else
 \def\next{\errmessage{Use \string\Bbb\space only in math
mode}}\fi\next}
\def\Bbb@#1{{\Bbb@@{#1}}}
\def\Bbb@@#1{\fam\msbfam#1}
\catcode`\@=12

 \catcode`\@=11
\font\twelveeuf=eufm10 scaled 1100
\font\teneuf=eufm10
\font\nineeuf=eufm7 scaled 1100
\newfam\euffam
\textfont\euffam=\twelveeuf  \scriptfont\euffam=\teneuf
  \scriptscriptfont\euffam=\nineeuf
\def\euf@{\hexnumber@\euffam}
\def\frak{\relax\ifmmode\let\next\frak@\else
 \def\next{\errmessage{Use \string\frak\space only in math
mode}}\fi\next}
\def\frak@#1{{\frak@@{#1}}}
\def\frak@@#1{\fam\euffam#1}
\catcode`\@=12
\font\emi= cmmi10 scaled 1700 
\font\eightmi=cmmi10

\title{Axiom A maps are dense in the space\\  of unimodal maps
  in the {\emi C}\raise6pt\hbox{\eightmi k} topology}
\shorttitle{Axiom A maps}  
 \author{O.\ S.\ Kozlovski}
 \institutions{The Warwick University, Coventry, United Kingdom\\
{\eightpoint {\it E-mail address\/}: oleg@maths.warwick.ac.uk}}

\newcommand{\be}{\begin} 
\newcommand{\en}{\end}

\newcommand{\Bg}{{\frak B}}

\newcommand{\ie}{{i.e.\ }}

\newcommand{\sto}[1]{{\displaystyle\mathop{#1}^{\circ}}}

\newcommand{\bpr}{ $\lhd$ } 
\newcommand{\epr}{ {\hfill~$\rhd$}  \vspace{1.5mm}}
 \newenvironment{plm}{\bpr}{\epr}

\newcommand{\comment}[1]{}


\newcommand{\la}{\lambda} \newcommand{\sm}{\setminus}
\newcommand{\Oi}{\mathop{{\it O}}}
\newcommand{\oi}{\mathop{{\it o}}}

\newcommand{\C}{\mathop{{{{\rm a}}}}\nolimits} 
\newcommand{\A}{\mathop{{{{\rm A}}}}\nolimits} 
\newcommand{\D}{\mathop{{{{\rm b}}}}\nolimits}
\newcommand{\B}{\mathop{{{{\rm B}}}}\nolimits}

\newcommand{\notsubset}{{\not \subset}}

\newcommand{\Cp}{{{{{\Bbb C}}}}}
\newcommand{\Rp}{{{{{\Bbb R}}}}}

\newcommand{\rc}{{\cal P}}


\newcommand{\eps}{\varepsilon}

\centerline{\bf Abstract}
\vglue12pt
  In this paper we   prove $C^k$ structural stability conjecture for
  unimodal maps. In other words, we shall prove that Axiom A maps are
  dense in the space of $C^k$ unimodal maps in the $C^k$ topology.
  Here $k$ can be $1,2,\ldots,\infty,\omega$.

 \section{Introduction}

1.1. {\it The structural stability conjecture}. The structural stability conjecture was and remains one of the
most interesting and important open problems in the theory of dynamical
systems.  This conjecture states that a dynamical system is
structurally stable if and only if it satisfies Axiom A and the
transversality condition. In this paper we prove this conjecture in
the simplest nontrivial case, in the case of smooth unimodal maps.
These are maps of an interval with just one critical turning point.

To be more specific let us recall the definition of Axiom A maps:

\numbereddemo{Definition}
  Let $X$ be an interval. We say that a $C^k$ map $f: X \hookleftarrow$
  satisfies the {\it Axiom} A conditions if:
  \begin{itemize}
  \item $f$ has finitely many hyperbolic periodic attractors,
  \item the set $\Sigma(f)= X \setminus \Bg(f)$ is hyperbolic,
    where $\Bg(f)$ is a union of the basins of attracting periodic
    points.
  \end{itemize}

\enddemo

This is more or less a classical definition of the Axiom A maps;
however in the case of $C^2$ one-dimensional maps Ma\~{n}\`e has
proved that a $C^2$ map satisfies Axiom A if and only if all its
periodic points are hyperbolic and the forward iterates of all its
critical points converge to some periodic attracting points.

It was proved many years ago that Axiom A maps are $C^2$ structurally
stable if the critical points are nondegenerate and the ``no-cycle''
condition is fulfilled (see, for example, \cite{dMvS}). However the
opposite question ``Does structural stability imply Axiom A?''
appeared to be much harder.  It was conjectured that the answer to
this question is affirmative and it was assigned the name
``structural stability conjecture''. So, the main result of this paper
is the following theorem:

\specialnumber{A}
\proclaim{Theorem} \label{AA}
  Axiom {\rm A} maps are dense in the space of $C^\omega(\Delta)$ unimodal
  maps in the $C^\omega(\Delta)$ topology {\rm (}$\Delta$ is an arbitrary
  positive number\/{\rm ).}\/
\endproclaim

Here $C^\omega(\Delta)$ denotes the space of real analytic functions
defined on the interval which can be holomorphically extended to a
$\Delta$-neighborhood of this interval in the complex plane.

Of course, since analytic maps are dense in the space of smooth maps
it immediately follows that $C^k$ unimodal Axiom A  maps  are dense in
the space of all unimodal maps in the $C^k$ topology, where
$k=1,2,\ldots,\infty$.

This theorem, together with the previously mentioned theorem, clearly
implies the structural stability conjecture:

\specialnumber{B}
\proclaim{Theorem}
  A $C^k$ unimodal map $f$ is $C^k$ structurally stable if and only if
  the map $f$ satisfies the Axiom {\rm A} conditions and  its critical point is
  nondegenerate and nonperiodic{\rm ,}
  $k=2,\ldots,\infty,\omega$.\footnote{If $k=\omega$, then one should
    consider the space $C^\omega(\Delta).$}
\endproclaim

Here the critical point is called {\it nondegenerate} if the second
derivative at the point is not zero.

In this theorem the number $k$ is greater than one because any
unimodal map can be $C^1$ perturbed to a nonunimodal map and, hence,
there are no $C^1$ structurally stable unimodal maps (the topological
conjugacy preserves the number of turning points). For the same
reason the critical point of a structurally stable map should be
nondegenerate.

In fact, we will develop tools and techniques which give more detailed
results.  In order to formulate them, we need the following definition: The map $f$ is {\it regular} if either
the $\omega$-limit set of its critical point $c$ does not contain neutral periodic points
or the $\omega$-limit set of $c$ coincides with the orbit of some
neutral periodic point. For example, if the map has negative
Schwarzian derivative, then this map is regular. Regular maps are
dense in the space of all maps (see Lemma~\ref{lm:regular}). We will
also show that if the analytic map $f$ does not have neutral periodic
points, then this map can be included in a family of regular analytic
maps.

\specialnumber{C}
\proclaim{Theorem} \label{real_anal} \hskip-8pt
  Let $X$ be an interval and $f_\la: X \hookleftarrow$ be an analytic
  family of analytic unimodal regular maps with a nondegenerate
  critical point{\rm ,}\break $\la \in \Omega \subset \Rp^N$ where $\Omega$ is a
  open set. If the family $f_\lambda$ is nontrivial in the sense that
  there exist two maps in this family which are not combinatorially
  equivalent{\rm ,} then Axiom {\rm A} maps are dense in this family.  Moreover{\rm ,}
  let $\Upsilon_{\la_0}$ be a subset of $\Omega$ such that the maps
  $f_{\la_0}$ and $f_{\la'}$ are combinatorially equivalent for $\la'
  \in \Upsilon_{\la_0}$ and the iterates of the critical point of
  $f_{\lambda_0}$ do not converge to some periodic attractor. Then the
  set $\Upsilon_{\la_0}$ is an analytic variety. If $N=1${\rm ,} then
  $\Upsilon_{\lambda_0}\cap Y${\rm ,} where the closure of the interval $Y$
  is contained in $\Omega${\rm ,} has finitely many connected components.
\endproclaim

Here we say that two unimodal maps $f$ and $\hat f$ are
{\it combinatorially equivalent} if there exists an order-preserving
bijection   $h:\cup_{n\geq 0}f^n(c)\to \cup_{v\geq 0}\hat f(\hat c) $
such that $h(f^n(c))=\hat f^n(\hat c)$ for all $n\geq 0$, where $c$ and
$\hat c$ are critical points of $f$ and $\hat f$. In the other words,
$f$ and $\hat f$ are combinatorially equivalent if the order of their
forward critical orbit is the same. Obviously, if two maps are
topologically conjugate, then they are combinatorially equivalent.


Theorem~\ref{AA} gives only global perturbations of a given map.
However, one can want to perturb a map in a small neighborhood of a
particular point and to obtain a nonconjugate map. This is also
possible to do and will be considered in a forthcoming paper. (In fact, all the tools and strategy of the proof
will be the same as in this paper.)

\demo{{\rm 1.2.} Acknowledgments}
First and foremost, I would like to thank\break S.~van~Strien for his
helpful suggestions, advice and encouragement. Special thanks go to
W.~de~Melo who pointed out that the case of maps having neutral
periodic points should be treated separately. His constant feedback
helped to improve and clarify the presentation of the paper.

G.~{\' S}wi{\c{a}}tek explained to me results on the quadratic family and our
many discussions clarified many of the concepts used here.
J.~Graczyk, G.~Levin and M.~Tsuji gave me helpful feedback at talks
that I gave during the International Congress on Dynamical Systems at
IMPA in Rio~de~Janeiro in 1997 and during the school on dynamical systems
in Toyama, Japan in 1998. I also would like to thank D.V.~Anosov,
M.~Lyubich, D.~Sands and E.~Vargas for their useful comments.

This work has been supported by the Netherlands Organization for
Scientific Research (NWO).
\enddemo

1.3. {\it Historical remarks}.  
The problem of the description of the structurally stable dynamical
systems goes back to Poincar\'e, Fatou, Andronov and Pontrjagin.
The explicit definition of a structurally stable dynamical system was
first given by Andronov although he assumed one extra condition: the
$C^0$ norm of the conjugating homeomorphism had to tend to 0 when
$\epsilon$ goes to 0.

 Jakobson proved that Axiom A maps are dense in the $C^1$ topology,
\cite{Jak}. The $C^2$ case is much harder and only some partial
results are known.  Blokh and Misiurewicz proved that any map
satisfying the Collect-Eckmann conditions can be $C^2$ perturbed to an
Axiom A map, \cite{BM2}. In \cite{BM1} they extend this result to a
larger class of maps. However, this class does not include the
infinitely renormalizable maps, and it does not cover
nonrenormalizable maps completely.

Much more is known about one special family of unimodal maps:
quadratic maps $Q_c:x\mapsto x^2+c$. It was noticed by Sullivan that
if one can prove that if two quadratic maps $Q_{c_1}$ and $Q_{c_2}$
are topologically conjugate, then these maps are quasiconformally
conjugate, then this would imply that Axiom A maps are dense in the
family $Q$. Now this conjecture is completely proved in the case of
real $c$ and many people made contributions to its solution: Yoccoz
proved it in the case of the finitely renormalizable quadratic maps,
\cite{Yoc}; Sullivan, in the case of the infinitely renormalizable
unimodal maps of ``bounded combinatorial type'', \cite{Sul1},
\cite{Sul2}.  Finally, in 1992 there appeared a preprint by {\'
  S}wi{\c{a}}tek where this conjecture was shown for all real
quadratic maps. Later this preprint was transformed into a joint paper
with Graczyk \cite{GS}. In the preprint \cite{Lyu2} this result was
proved for a class of quadratic maps which included the real case as
well as some nonreal quadratic maps; see also \cite{Lyu4}. Another
proof was recently announced in \cite{Shi}.  Thus, the following
important rigidity theorem was proved:

\nonumproclaim{Theorem {\rm (Rigidity Theorem)}}
  If two quadratic non Axiom {\rm A} maps $Q_{c_1}$ and $Q_{c_2}$ are
  topologically conjugate ($c_1,c_2 \in \Rp$){\rm ,} then $c_1=c_2$.
\endproclaim

1.4. {\it Strategy of the proof}. Thus, we know that we can always perturb a quadratic map and change
its topological type if it is not an Axiom A map. We want to do the
same with  an arbitrary unimodal map of an interval. So the first
reasonable question one may ask is ``What makes  quadratic maps so
special''? Here is a list of major properties of the
quadratic maps which the ordinary unimodal maps do not enjoy:
\begin{itemize}
\item  Quadratic maps are analytic and they have nondegenerate
  critical point;
\item  Quadratic maps have negative Schwarzian derivative;
\item Inverse branches of quadratic maps have ``nice'' extensions to
  the complex plane (in terminology which we will introduce later we
  will say that the quadratic maps belong to the Epstein class);
\item Quadratic maps are polynomial-like maps;
\item The quadratic family is rigid in the sense that a quasiconformal
  conjugacy between two non Axiom A maps from this family implies that
  these maps coincide;
\item  Quadratic maps are regular.
\end{itemize}

We will have to compensate for the lack of these properties somehow. 

First, we   notice that since the analytic maps are dense in the
space of $C^k$ maps it is sufficient to prove the $C^k$ structural
stability conjecture only for analytic maps, i.e.,  when $k$ is
$\omega$.  Moreover, by the same reasoning we can assume that the
critical point of a map we want to perturb is nondegenerate.

The negative Schwarzian derivative condition is a much more subtle
property and it provides the most powerful tool in one-dimensional
dynamics. There are many theorems which are proved only for maps with
negative Schwarzian derivative. However, the tools described in
\cite{Koz} allow us to forget about this condition! In fact, any
theorem proved for maps with negative Schwarzian derivative can be
transformed (maybe, with some modifications) in such a way that it is
not required that the map have negative Schwarzian derivative anymore. Instead of the negative Schwarzian
derivative the map will have to have a nonflat critical point.

In the first versions of this paper, to get around the Epstein class,
we needed to estimate the sum of lengths of intervals from an orbit of
some interval. This sum is small if the last interval in the orbit is
small. However, Lemma~2.4 in \cite{dFdM} allows us to estimate the shape
of pullbacks of disks if one knows an estimate on the sum of lengths
of intervals in some power greater than 1. Usually such an estimate is
fairly easy to arrive at and in the present version of the paper we do not
need estimates on the sum of lengths any more.

Next, the renormalization theorem will be proved; i.e.\  we
will prove that for a given unimodal analytical map with a
nondegenerate critical point there is an induced holomorphic
polynomial-like map, Theorem~\ref{p_like}.  For infinitely
renormalizable maps this theorem was   proved in \cite{LvS}.
For finitely renormalizable maps we will have to generalize the notion
of polynomial-like maps, because one can show that the classical
definition does not work in this case for all maps.

Finally, using the method of quasiconformal deformations, we will
construct a perturbation of any given analytic regular map and show
that any analytic map can be included in a nontrivial analytic family
of unimodal regular maps.

If the critical point of the unimodal map is not recurrent, then
either its forward iterates converge to a periodic attractor (and if
all periodic points are hyperbolic, the map satisfies Axiom A) or this
map is a so-called Misiurewicz map. Since in the former case we have
nothing to do the only interesting case is the latter one. However,
the Misiurewicz maps are fairly well understood and this case is
really much simpler than the case of maps with a recurrent critical
point. So, usually we will concentrate on the latter, though  the case of Misiurewicz maps is also
considered. 
 
We have tried to keep the exposition in such a way that all section of
the paper are as independent as possible. Thus, if the reader is
interested only in the proofs of the main theorems, believes that maps
can be renormalized as described in Theorem~\ref{p_like} and is
familiar with standard definitions and notions used in one-dimensional
dynamics, then he/she can start reading the paper from
Section~\ref{ch:Ck}.

\demo{{\rm 1.5.} Cross-ratio estimates}
Here we briefly summarize some known facts about cross-ratios which we
will use intensively  throughout the paper.

There are several types of cross-ratios which work more or less in the
same way. We will use just a standard cross-ratio which is given by
the formula:
$$
\D(T,J)=\frac{|J||T|}{|T^-||T^+|}
$$
where $J\subset T$ are intervals and $T^-$, $T^+$ are connected
components of $T\setminus J$.

Another useful cross-ratio (which is in some sense degenerate) is the
following:
$$
\C(T,J)=\frac{|J||T|}{|T^-\cup J||J\cup T^+|}
$$
where the intervals $T^-$ and $T^+$ are defined as before.

If $f$ is a map of an interval, we will measure how this map distorts
the cross-ratios and introduce the following notation:
\begin{eqnarray*}
  \B(f,T,J)&=& \frac{\D(f(T),f(J))}{\D(T,J)}\\[6pt]
  \A(f,T,J)&=& \frac{\C(f(T),f(J))}{\C(T,J)}.
\end{eqnarray*}

It is well-known that maps having negative Schwarzian derivative
increase the cross-ratios: $ \B(f,T,J)\geq 1$ and $ \A(f,T,J)\geq 1$ if
$J\subset T$, $f|_T$ is a diffeomorphism and the $C^3$ map $f$ has
negative Schwarzian derivative. It turns out that if the map $f$ does
not have negative Schwarzian derivative, then we also have an estimate
on the cross-ratios provided the interval $T$ is small enough. This
estimate is given by the following theorems (see \cite{Koz}):
\enddemo

\advance\theoremcount by -1
\proclaim{Theorem} \label{cr1} 
  Let $f: X \hookleftarrow$ be a $C^3$ unimodal map of an interval to
  itself with a nonflat nonperiodic critical point and suppose that
  the map $f$ does not have any neutral periodic points. Then there
  exists a constant $C_1>0$ such that if $M$ and $I$ are
  intervals{\rm ,} $I$ is a subinterval of $M${\rm ,} $f^n|_M$ is monotone and
  $f^n(M)$ does not intersect the immediate basins of periodic
  attractors{\rm ,} then
  \begin{eqnarray*}
    \A(f^n,M,I)&>&\exp(-C_1\, |f^n(M)|^2),\\[6pt]
    \B(f^n,M,I)&>&\exp(-C_1\, |f^n(M)|^2).
  \end{eqnarray*} 
\endproclaim 
\pagebreak

Fortunately, we will usually deal only with maps which have no
neutral periodic points because such maps are dense in the space of
all unimodal maps. However, at the end we will need some estimates for
maps which do have neutral periodic points and then we will use
another theorem (\cite{Koz}):

\proclaim{Theorem} \label{nsch} 
  Let $f: X \hookleftarrow$ be a $C^3$ unimodal map of an interval to
  itself with a nonflat nonperiodic critical point. Then there
  exists a nice\footnote{The definition of nice intervals is given in
    the next subsection.} interval $T$ such that the first entry map
  to the interval $f(T)$ has negative Schwarzian derivative.
\endproclaim

1.6. {\it Nice intervals and first entry maps}.
In this section we introduce some definitions and notation.

The {\it basin} of a periodic attracting orbit is a set of points
whose iterates converge to this periodic attracting orbit. Here the
periodic attracting orbit can be neutral and it can attract points
just from one side. The {\it immediate} basin of a periodic attractor
is a union of connected components of its basin whose contain points
of this periodic attracting orbit. The union of immediate basins of
all periodic attracting points will be called the {\it immediate
  basin of attraction} and will be denoted by $\Bg_0$.

We say that the point $x'$ is {\it symmetric} to the point $x$ if
$f(x)=f(x')$. In this case we call the interval $[x,x']$
{\it symmetric} as well.  A symmetric interval $I$ around a critical
point of the map $f$ is called {\it nice} if the boundary points of
this interval do not return into the interior of this interval under
iterates of $f$. It is easy to check that there are nice intervals of
arbitrarily small length if the critical point is not periodic.

Let $T\subset X $ be a nice interval and $f: X \hookleftarrow$ be a
unimodal map. $R_T:U\to T$ denotes the first entry map to the interval
$T$, where the open set $U$ consists of points which occasionally
enter the interval $T$ under iterates of $f$. If we want to consider
the first return map instead of the first entry map, we will write
$R_T|_T$.  If a connected component $J$ of the set $U$ does not contain the
critical point of $f$, then $R_T:J\to T$ is a diffeomorphism of the
interval $J$ onto the interval $T$. A connected component of the set
$U$ will be called a {\it domain} of the first entry map $R_T$, or a
{\it domain} of the nice interval $T$. If $J$ is a domain of $R_T$,
the map $R_T:J\to T$ is called a {\it branch} of $R_T$. If a domain
contains the critical point, it is called {\it central}.

Let $T_0$ be a small nice interval around the critical point $c$ of
the map $f$. Consider the first entry map $R_{T_0}$ and its central
domain.  Denote this central domain as $T_1$. Now we can consider the
first entry map $R_{T_1}$ to $T_1$ and denote its central domain as
$T_2$ and so on. Thus, we get a sequence of intervals $\{T_k\}$ and a
sequence of the first entry maps $\{R_{T_k}\}$.

We will distinguish several cases. If $c \in R_{T_k}(T_{k+1})$, then
$R_{T_k}$ is called a {\it high} return and if $c \notin
R_{T_k}(T_{k+1})$, then $R_{T_k}$ is a {\it low} return. If
$R_{T_k}(c) \in T_{k+1}$, then $R_{T_k}$ is a {\it central} return and
otherwise it is a {\it noncentral} return.

The sequence $T_0\supset T_1\supset \cdots$ can converge to some
nondegenerate interval~$\tilde T$. Then the first return map
$R_{\tilde T}|_{\tilde T}$ is again a unimodal map which we call a
{\it renormalization} of $f$ and in this case the map $f$ is called
{\it renormalizable} and the interval $\tilde T$ is called a {\it
  restrictive} interval. If there are infinitely many intervals such
that the first return map of $f$ to any of these intervals is
unimodal, then the map $f$ is called {\it infinitely renormalizable}.

Suppose that $g: X \hookleftarrow$ is a $C^1$ map and suppose that
$g|_J:J\to T$ is a diffeomorphism of the interval $J$ onto the
interval $T$. If there is a larger interval $J'\supset J$ such that
$g|_{J'}$ is a diffeomorphism, then we will say that the range of the
map $g|_J$ can be {\it extended} to the interval $g(J')$.

We will see that any branch of the first entry map can be decomposed
as a quadratic map and a map with some definite extension.

\specialnumber{1.1}\proclaim{Lemma} \label{exten_centr} 
  Let $f$ be a unimodal map{\rm ,} $T$ be a nice interval{\rm ,} $J$ be its
  central domain and $V$ be a domain of the first entry map to $J$
  which is disjoint from $J${\rm ,} i.e.\  $V\cap J=\emptyset$. Then
  the range of the map $R_J:V\to J$ can be extended to $T$.
\endproclaim

This is a well-known lemma; see for example \cite{dMvS} or \cite{Koz}.

We say that an interval $T$ is a  $\tau$-scaled neighborhood of the
interval $J$, if $T$ contains $J$ and if each component of $T\setminus
J$ has at least length $\tau|J|$.

\vglue-6pt
\section{Decay of geometry}
\label{sec:decay}

\vglue-6pt

In this section we state an important theorem about the exponential
``decay of geometry''.  We will consider {\it unimodal
  nonrenormalizable} maps with a recurrent {\it quadratic} critical
point. It is known that in the multimodal case or in the case of a
degenerate critical point this theorem does not hold.

Consider a sequence of intervals $\{T_0,T_1,\ldots\}$ such that the
interval $T_0$ is nice 
and the interval $T_{k+1}$ is a central domain of the first entry map
$R_{T_k}$.  Let $\{k_l,\,l=0,1,\ldots\}$ be a sequence such that
$T_{k_l}$ is a central domain of a noncentral return. It is easy to
see that since the map $f$ is nonrenormalizable the sequence
$\{k_l\}$ is unbounded and the size of the interval $T_k$ tends to 0
if $k$ tends to infinity. 

The decay of the ratio $\frac{|T_{k_l+1}|}{|T_{k_l}|}$ will play an
important role in the next section.

\proclaim{Theorem} \label{exp_decay} \label{thr:decay}
  Let $f$ be an analytic unimodal nonrenormalizable map with a
  recurrent quadratic critical point and without neutral periodic
  points.  Then the ratio $\frac{|T_{k_l+1}|}{|T_{k_l}|}$ decays
  exponentially fast with $l$.
\endproclaim 

This result was suggested in \cite{Lyu3} and it has been proven in
\cite{GS} and \cite{Lyu4} in the case when the map is quadratic or
when it is a box mapping. To be precise we will give the statement of
this theorem below, but first we introduce the notion of a box
mapping.

\demo{Definition {\rm 2.1}}
  Let $A\subset \Cp$ be a simply connected Jordan domain,\break $B \subset
  A$ be a domain   each of whose connected components is a simply connected
  Jordan domain and let $g:B \to A$ be a holomorphic map. Then $g$ is
  called {\em a holomorphic box mapping} if the following assumptions
  are satisfied:
  \begin{itemize}
  \item $g$ maps the boundary of a connected component of $B$ onto the
    boundary of $A$,
  \item There is one component of $B$ (which we will call {\em a
      central domain}) which is mapped in the $2$-to-$1$ way onto the
    domain $A$ (so that there is a critical point of $g$ in the central
    domain),
  \item All other components of $B$ are mapped univalently onto $A$ by
    the map~$g$,
  \item The iterates of the critical point of $g$ never leave the domain
    $B$.
  \end{itemize}
\enddemo

In our case all holomorphic box mappings will be called {\it real} in the sense
that the domains $B$ and $A$ are symmetric with respect to the real
line and the restriction of $g$ onto the real line is real. 

We will say that a real holomorphic box mapping $F$ is {\it induced } by
an analytic unimodal map $f$ if any branch of $F$ has the form $f^n$.

We can repeat all constructions we used for a real unimodal map in the
beginning of this section for a real holomorphic box mapping.
Denote the central domain of the map $g$ as $A_1$ and consider the
first return map onto $A_1$. This map is again a real holomorphic box
mapping and we can again consider the first return map onto the domain
$A_2$ (which is a central domain of the first entry map onto $A_1$)
and so on. The definition of the central and noncentral returns and
the definition of the sequence $\{k_l\}$ can be literally transferred
to this case if $g$ is nonrenormalizable (this means that the
sequence $\{k_l\}$ is unbounded).

\proclaimtitle{\cite{GS}, \cite{Lyu4}}
\proclaim{Theorem}   \label{holom}
  Let $g:B \to A$ be a real holomorphic nonrenormalizable box mapping
  with a recurrent critical point and let the modulus of the annulus
  $A\sm \hat B$ be uniformly bounded from $0${\rm ,} where $\hat B$ is any
  connected component of the domain $B$. Then the ratio
  $\frac{|A_{k_l+1}|}{|A_{k_l}|}$ tends to $0$ exponentially fast{\rm ,} where
  $|A_k|$ is the length of the real trace of the domain $A_k$.
\endproclaim 

Here the real trace of the domain is just the intersection of this
domain with the real line.

So, if we can construct an induced box mapping, we will be able to
prove Theorem~\ref{exp_decay}. Fortunately, this construction has been
done in \cite{LvS} and in the less general case in \cite{GS}, \cite{Lyu3}.

\proclaim{Theorem}\label{thr:lvs}  
For any analytic unimodal map $f$ with a   nondegenerate critical point there exists an induced holomorphic
box
  mapping $F:B\to A$. Moreover{\rm ,} there exists a constant $C>0$ such
  that if $\hat B$ is a connected component of $B${\rm ,} then
  ${\rm mod}\,(A\setminus \hat B)>C$.
\endproclaim 

In fact, this theorem was proven in \cite{LvS} for infinitely
renormalizable maps in full generality and for the finitely
renormalizable maps satisfying two extra assumptions: $f$ has negative
Schwarzian derivative and $f$ belongs to the Epstein class (for
definition of the Epstein class see Appendix~5.2).
However, these conditions are not necessary any more. Indeed,
Theorem~\ref{thr:lvs} is a consequence of some estimates (usually
called ``complex bounds''). In \cite{LvS} these estimates are robust
in the following sense: if you change all constants involved by some
spoiling factor which is close to $1$, then the estimates still remain
true. Now, according to \cite{Koz} on small scales one has the
cross-ratio estimates as in the case of maps with negative Schwarzian
derivative, but with some spoiling factor close to $1$ (see
Theorems~\ref{cr1} and \ref{nsch}). Lemma~2.4 in \cite{dFdM} gives
estimates for the shape of pullbacks of disks and makes the Epstein
class condition superficial.  This lemma is formulated below in
Appendix~5.2 (Lemma~\ref{lm:fm}).  Thus, the combination
of Lemma~2.4 in \cite{dFdM}, the results of \cite{Koz} and of the
proof of the renormalization theorem in \cite{LvS} provides
Theorem~\ref{thr:lvs}. The outline of the proof is given in
Appendix~5.3. 

Theorem~\ref{exp_decay} is a trivial consequence of
Theorems~\ref{holom} and \ref{thr:lvs}.

\section{Polynomial-like maps} \label{ch:p_like}

The notion of polynomial-like maps was introduced by A.~Douady and\break
J.\ H.~Hubbard and was generalized several times after that. The main
advantage of using this notion is that one can work with a
polynomial-like map in the same way as if it was just a polynomial
map. We will use the following definition:

\numbereddemo{Definition}
  A holomorphic map $F:B \to A$ is called {\it polynomial-like} if it
  satisfies the following properties:
  \begin{itemize}
  \item $B$ and $A$ are domains in the complex plane, each having
    finitely many connected components; each connected component of
    $B$ or $A$ is a simply connected Jordan domain and $B$ is a subset
    of $A$.  The intersection of the boundaries of the domains $A$ and
    $B$ is empty or it is a forward invariant set which consists of
    finitely many points;
  \item The boundary of a connected component of $B$ is mapped onto the
    boundary of some connected component of $A$;
  \item There is one selected connected component $B^c$ of $B$ (which
    we will call {\it central}) such that the map $F|_{B^c}$ is
    $2$-to-$1$, and the central component $B^c$ is relatively compact in
    the domain $A$ (\ie $\bar{B^c}\subset A$);
  \item On the other connected components of $B$ the map $F$ is
    univalent.
  \end{itemize}
\enddemo

 \vglue-6pt

If the domains $A$ and $B$ are simply connected and the annulus
$A\setminus B$ is not degenerate, then a polynomial-like map
$F:B\to A$ is called a {\it quadratic-like} map.

We say that the polynomial-like map is {\em induced} by the unimodal
map $f$ if all connected components of the domains $A$ and $B$ are
symmetric with respect to the real line and the restriction of $F$ on
the real trace of any connected component of $B$ is an iterate of the
map $f$.

Notice a similarity between polynomial-like maps and holomorphic box
maps. There are two differences: in the case of the polynomial-like
map the domains $A$ and $B$ consist of several connected components
and in the case of the holomorphic box map the domain $A$ is simply
connected and the domain $B$ can consist of infinitely many connected
components. It is easy to see that if the critical point never leaves
$B$ under iterations of $F$, then the first return map of a
polynomial-like map to the connected component of $A$ which contains
the critical point is a holomorphic box map.

The main result of this section is that an analytic unimodal map can
be ``renormalized'' to obtain a polynomial-like map.

Before giving the statement of the theorem let us introduce the
following notation. $D_\phi(I)$ will denote a lens, \ie an
intersection of two disks of the same radius in such a way that two
points of the intersection of the boundaries of these disks are joined
by $I$ and the angle of this intersection at these points is $2\phi$.
See also Appendix~5.2 and Figure~1.
 \vglue-12pt
\figin{fig1}{1000}
\centerline{Figure 1. The lens $D_\phi(I)$} 

\advance\theoremcount by -1
\proclaim{Theorem} \label{p_like} 
  Let $f$ be an analytic{\rm ,} unimodal{\rm ,} not infinitely renormalizable map
  with a quadratic recurrent critical point and without neutral
  periodic points.  Then for any $\epsilon>0$ there exists a
  polynomial\/{\rm -}\/like map $F:B\to A$ induced by the map $f${\rm ,} and
  satisfying the following properties\/{\rm :}
  \begin{itemize}
  \item The forward orbit of the critical point under iterations of
    $F$ is contained in $B${\rm ;}
  \item $A$ is a union of finitely many lenses of the form
    $D_\phi(I)${\rm ,} where $I$ is an interval on the real line{\rm ,}
    $|I|<\epsilon$ and $0<\phi<\pi/4${\rm ;}
  \item If $F(x)\in A^c${\rm ,} then $B^x$ is compactly contained in $A^x${\rm ,}
    where $B^x$ and $A^x$ denote connected components of $B$ and $A$
    containing $x$ and $A^c$ denotes a connected component of $A$
    containing the critical point $c$ {\rm (}\/\ie $\bar{B^x}\subset A^x${\rm ,}
    where $\bar{B^x}$ is the closure of $B^x${\rm );}
  \item Boundaries of connected components of $B$ are piecewise smooth
    curves{\rm ;}
  \item If $a\in \partial A \cap \partial B${\rm ,} then the boundaries of
    $A$ and $B$ at $a$ are not smooth{\rm ;} however if we consider a smooth
    piece of the boundary of $A$ containing $a$ and the corresponding
    smooth piece of the boundary of $B${\rm ,} then these pieces have the
    second order of tangency {\rm (}\/see Figure~{\rm 2);}
   \item If $B^{x_1}\cap B^{x_2}=\emptyset$ and $b\in \partial
    B^{x_1}\cap \partial B^{x_2}${\rm ,} then the boundaries of $B^{x_1}$
    and $B^{x_2}$ are not smooth at the point $b$ and not tangent to
    each other\/{\rm ; }\/
  \item For any $x\in B${\rm ,}
    $$
    \frac{|B^x|}{|A^x|}<\epsilon,
    $$
    where $|B^x|$ denotes the length of the real trace of $B^x${\rm ;}
  \item If $x\in B$ and $F|_{B^x}=f^n${\rm ,} then $f^i(x)\notin A^c$ for
    $i=1,\dots,n-1${\rm ;}
  \item $f(c)\notin A${\rm ;}
  \item When $a\in \partial A$ is a point closest to the critical value
    $f(c)${\rm ,} then 
    $$\frac{|f(B^c)|}{|a-f(c)|}<\epsilon.$$
  \end{itemize}
\endproclaim 
  \vglue-36pt
\figin{fig2}{900}
\begin{quote}
{Figure 2. {A fragment of the domain of definition of a
  polynomial-like map}}
\end{quote}

If the map $f$ is infinitely renormalizable, we will use a  much simpler
statement. 

\proclaimtitle{\cite{LvS}}
\proclaim{Theorem} \label{thr:ren_minimal}
  Let $f$ be an analytic unimodal infinitely renormalizable map with a
  quadratic critical point. Then there exists a quadratic\/{\rm -}\/like map
  $F:B\to A$ induced by $f$ such that the forward orbit of $c$ under
  iterates of $F$ is contained in $B$.
\endproclaim 

The proof of Theorem~3.1 will occupy the rest of this
section.

\demo{{\rm 3.1.} The real and complex bounds}
In this subsection we give two technical lemmas. 
\enddemo

\specialnumber{3.1}\proclaim{Lemma} \label{real_b}
  Let $f$ be a $C^3$ nonrenormalizable unimodal map with a quadratic
  recurrent critical point. Then for any $\epsilon>0$ there exists
  $\delta>0$ such that if $T_0$ is a sufficiently small nice interval{\rm ,}
  $T_1$ is a central domain of $T_0${\rm ,} $T_2$ is a central domain of
  $T_1$ and $\frac{|T_1|}{|T_0|}<\delta${\rm ,} then the following holds\/{\rm :} When
  $T_1'$ is a domain of $R_{T_1}$ containing the critical value $f(c)$
  {\rm (}\/see Fig.\ {\rm 3),} then
  $$
  \frac{|T_1'|}{|f(T_1)|} < \epsilon.
  $$
\endproclaim

\centerline{\BoxedEPSF{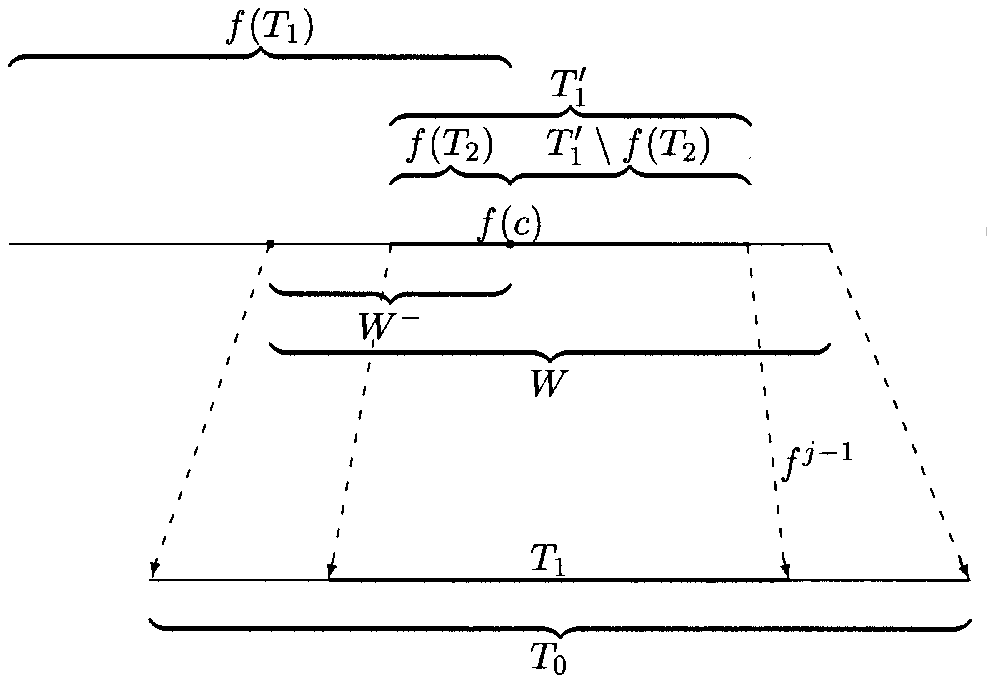 scaled 800}}
\vglue-24pt
\centerline{Figure 3. The map $f^{j-1}$.}
\vglue12pt

$\lhd$  Let $R_{T_1}|_{T_2}=f^j$. The range of the map $f^{j-1}:T_1'\to T_1$
can be extended to the interval $T_0$ (Lemma~\ref{exten_centr});
i.e.,  there is an interval $W$ such that\break $f^{j-1}:W\to T_0$ is
a diffeomorphism, $T_1'\subset W$ and $f^{j-1}(W)=T_0$. Denote the
components of $W\sm (T_1'\sm f(T_2))$ as $W^-$ and $W^+$ in such a way
that the interval $f(T_2)$ is a subset of the interval $W^-$. It is
easy to see that the interval $f(T_1)$ contains the interval $W^-$.
Applying Theorem~\ref{cr1} we obtain the following bounds:
\begin{eqnarray*}
  \frac{|T_1'|}{|f(T_1)|}
  &\leq&
  \frac{|T_1'|}{|W^-|} \leq 
  \D(W,T_1')\\
  &\leq&
  \D(T_0,T_1)  \leq 
 C_2 \frac{4\delta}{(1+\delta)^2}
\end{eqnarray*}
where the constant $C_2$ is close to 1 if the interval $T_0$ is
sufficiently small.   


\specialnumber{3.2}\proclaim{Lemma} \label{lm:complex2}
  Let $f$ be an analytic unimodal map.  For any $\phi_0\in (0,\pi)$
  and $K>0$ there are constants $\phi\in (0,\phi_0)$ and
  $C_3>0$ such that if $f^n|_V$ is monotone{\rm ,}
  $|f^i(V)|<C_3$ for $i=0,\dots,n$ and $\sum_{i=0}^n|f^i(V)|<K${\rm ,}
  then
  $$f^{-n}(D_\phi(f^n(V)))\subset  D_{\phi_0}(V).$$
\endproclaim

This lemma is a simple consequence of Lemma~5.2 in \cite[p.~487]{dMvS}. One can also use Lemma~2.4 in
\cite{dFdM} which gives better estimates (see Lemma~\ref{lm:fm}).

\demo{{\rm 3.2.} Construction of the induced polynomial\/{\rm -}\/like
  map}
\enddemo

{\it Proof of Theorem} 3.1. If the $\omega$-limit set of the critical point is
minimal (we say that the forward invariant set is {\it minimal} if it
closed and has no proper closed invariant subsets), then one can
construct the polynomial-like map in a much simpler way than is
given here.  In fact, it is a consequence of Theorem~2.3.
For example, the domain $A$ in this case is simply connected. However,
if the $\omega$-limit set of the critical point contains intervals,
the domain $A$ cannot be connected if we want the domain $B$ to
contain finitely many connected components.

Letting $\phi_0=\pi/4$, $K=|X|$, we  apply Lemma~\ref{lm:complex2} to
the map $f$ and obtain two constants $\phi$ and $C_3$. 

On the other hand, for this constant $\phi$ there is a constant
$\tau_1$ such that if an interval $J$ contains a
$\tau_1$-scaled neighborhood of an interval $I$, then
$D_{\pi/4}(I)\subset D_\phi(J)$.

Take a nice interval $T_0$ such that 
\begin{itemize}
\item $|T_0|<\epsilon$;
\item The boundary points of $T_0$ are eventually mapped by $f$ onto
  some repelling periodic point and $T_0$ is disjoint from the
  immediate basin of attraction $\Bg_0$;
\item The central domain $T_1$ of $T_0$ is so small that
  $\frac{|T_1'\setminus f(T_2)|}{|f(T_1)|}<\min(\frac 12 \tan^2 \frac
  \phi 2,\epsilon)$, where $T_2$ is a central domain of $T_1$ and
  $T_1'$ is a domain of $R_{T_1}$ containing the critical value (due
  to Theorem~\ref{thr:decay} the ratio $\frac{|T_1|}{|T_0|}$ can be
  made arbitrarily small and then we can apply Lemma~\ref{real_b});
\item If $f^n|_V$ is monotone and $f^n(V)\subset T_1$, then
  $|V|<C_3$ (the existence of such an interval $T_0$ follows from
  the absence of wandering intervals, for details see Lemma~5.2 in
  \cite{Koz});
\item Moreover, the ratio $\frac{|T_1|}{|T_0|}$ should be so small
  that if $f^n|_V$ is monotone and $f^n(V)=T_0$, then $V$ contains a
  $\tau_1$-scaled neighborhood of the pullback $f^{-n}(T_1)$ and
  $\frac{|f^{-n}(T_1)|}{|V|}<\epsilon$ (indeed, if
  $\frac{|T_1|}{|T_0|}$ is small, then the cross-ratio $\D(T_0,T_1)$
  is also small, the pullback can only slightly increase this cross-ratio, so that  $\D(V,f^{-k}(T_1))$ is
small; hence $f^{-k}(T_1)$ is deep
  inside $V$).
\end{itemize}

Let $\Bg_0$ be the immediate basin of attraction.  It is known that
the periods of attracting or neutral periodic points are bounded
(\cite{MdMvS}). Hence, the set $ X  \sm \bar \Bg_0$ consists of
finitely many intervals (as usual $\bar \Bg_0$ is a closure of
$\Bg_0$).  Some points of the interval $ X $ are mapped to the
immediate basin of attraction after some iterates of $f$. Obviously,
for a given $n$, the set $\{x\in  X :\, f^n(x)\notin \bar \Bg_0\}$
consists of finitely many intervals as well.

Just to fix the situation let us suppose that the map
$f:X\hookleftarrow$ first increases and then decreases.  Let
$P_n=\{x\in (\partial_- X,f(\partial T_1)):\, f^i(x)\notin \bar
T_1\cup \bar \Bg_0 \hbox{ for } i=0,\ldots,n\}$, where $\partial_- X$
denotes the left boundary point of $X$. The set $P_n$ consists of
finitely many intervals and the lengths of these intervals tend to
zero as $n\to \infty$ (otherwise we would have a wandering interval).
All the boundary points of $P_n$ are eventually mapped onto some
periodic points.  Moreover, the set of these periodic points is finite
and does not depend on $n$. Denote the union of this set and
$\omega(\partial T_1)$ (which is an orbit of a periodic point by the
choice of $T_0$) by $E$.  Let $a\in E$ be a periodic orbit of
period $k$. Then there exists a neighborhood of $a$ where the map
$f^k$ is holomorphically conjugate to a linear map.  This implies that
if $V$ is a sufficiently small interval and $a$ is its boundary point,
then $f^{-2k}(D_\phi(V))\subset D_\phi(V)$; hence
$f^{-2k(i+1)}(D_\phi(V))\subset f^{-2ki}(D_\phi(V))$ for
$i=0,1,\ldots$ and the size of $f^{-2ki}(D_\phi(V))$ tends to zero.

Due to a theorem of M\~{a}n\`e there exist two constants $C_4>0$
and $\tau_2>1$ such that if $x\in P_n$, then
$Df^i(x)>C_4\tau_2^i$ for $i=0,\ldots,n$ (see Theorem~5.1 in
\cite[p.\ 248]{dMvS}). Therefore there exists a constant
$C_5>0$ such that if $V\subset P_n$ is an interval, and
$|f^n(V)|<C_5$, then $|f^i(V)|<C_3$ for $i=0,\ldots,n$, and
$\sum_{i=0}^n |f^i(V)|<|X|$.

Let $m$ be so large that if $V$ is a connected component of $P_m$,
then $|V|<\min(C_5,\epsilon)$ and, moreover, if $V$ contains a
periodic point in its boundary, then $V$ is so small that the lens
$D_\phi(V)$ satisfies the properties described above (so it should be
in a neighborhood of this periodic point where the map can be
linearized and the size of the pullback of $D_\phi(V)$ along this
periodic orbit tends to zero).

Once we have fixed the integer $m$, we are not going to change it and
thus we will suppress the dependence of $P_m$ on $m$.

Let $S$ be a union of the boundary of the set $P$ and the forward
orbit of $\partial T_1$. Notice that $S$ is a finite forward invariant
set. The partition of the set $P\cup T_1$ by points of $S$ we denote
by $\rc$. Finally, let $A=\bigcup_{V\in\rc} D_\phi(V)$. The set $A$
will be the range of the polynomial-like map we are constructing.

Let $\Sigma$ be a closure of all points on the real line whose
$\omega$-limit set contains the critical point.
For any point $x\in\Sigma'=\Sigma\cap (\bar P\cup \bar T_1)$ such that
$f^i(x)\notin E$ for any $i>0$, we will construct an interval $I(x)$
and an integer $n(x)$ such that $x\in I(x)$, $f^{n(x)}(I(x))\in\rc$
and $f^{-n(x)}(D_\phi(\rc(f^{n(x)}(x))))\in D_\phi(\rc(x))$, where
$\rc(x)$ denotes an element of the partition containing the point $x$.
If the point $x\in \Sigma'$ is eventually mapped to some point of
$E$ and on both sides of $x$ there are points of $\Sigma'$
arbitrarily close to $x$, then we will construct two intervals
$I_-(x)$ and $I_+(x)$ on both sides of $x$ and two integers $n_-(x)$
and $n_+(x)$ with similar properties. If $f^i(x)\in E$ but there are no points of $\Sigma'$ on one side of
$x$ close to $x$, only intervals on the side containing points of $\Sigma'$ will be constructed.
Finally, if $x\in T_2$, we will put $I(x)=T_2$ and $n(x)$ will be a
minimal positive integer such that $f^{n(x)}(x)\in T_1$. In
this case $f^{n(x)}(I(x))\subsetneq T_1$ and so $f^{n(x)}(I(x))\notin
\rc$, however as we will see below $f^{-n(x)}(D_\phi(T_1))\subset
D_\phi(T_1)$.


First, we are going to construct these intervals and integers for a point
$x$ whose orbit contains points of the set $S$, where $S$ is a set of
boundary points of $P$. In this case some iterate of $x$ lands on a
periodic point $a\in E$; i.e., $f^k(x)=a\in E$. For   simplicity
let us assume that $a$ is just a fixed point and that its multiplier
is positive. Let $J$ be an interval of $\rc$ containing $a$ (there are
at most two such intervals). Because of the choice of $m$
\comment{(see item~\ref{n1})} we know that
$f|_J^{-1}(D_\phi(J))\subset D_\phi(J)$ and since $D_\phi(J)$ is in
the neighborhood of $a$ where the map $f$ can be linearized, the sizes
of domains $f|_J^{-i}(D_\phi(J))$ shrink to zero when $i\to +\infty$.
Thus, there exists $i_0$ such that 
$$f^{-k}\circ
f|_J^{-i_0}(D_\phi(J))\subset D_\phi(J')$$ 
and
$$\frac{|f^{-k}\circ f|_J^{-i_0}(J)|}{|J'|}<\epsilon,$$
where $J'$ is
just $\rc(x)$ if $x\notin S$ and $J'$ is one of the intervals of $\rc$
which contains $x$ on its boundary if $x\in S$.  We put
$I_-(x)=f^{-k}\circ f|_J^{-i_0}(J)$ and $n_-(x)=k+i_0$. If there is another interval from $\rc$ containing $a$
in its boundary, we can repeat the procedure and get the interval $I_+(x)$ and the integer
$n_+(x)$; otherwise we are finished in this case.

Now let us consider the case when $f^i(x)\notin S$ for all $i>0$. This
case we divide in several subcases.
\vglue6pt
If $x\in T_2$, then     $I(x)=T_2$ and $n(x)$ is a minimal
positive integer such that $f^{n(x)}(T_2)\subset T_1$; i.e.,
$R_{T_1}|_{T_2}=f^{n(x)}$. Let $T_1'$ be an interval around the
critical value $f(c)$ such that $f^{n(x)-1}(T_1')=T_1$ (see
Figure 3).  The pullback of a lens $D_\phi(T_1)$ by
$f^{-(n(x)-1)}$ is contained in $D_{\pi/4}(T_1')$ (indeed, by the
choice of $T_0$ we know that all intervals in the orbit
$\{f^i(T_1'),\, i=0,\dots,n(x)\}$ are small and they are disjoint; so we
can apply Lemma~\ref{lm:complex2}).  Near the critical point the map
$f$ is almost quadratic (if $T_0$ is small enough) and because of the
choice of $T_0$ the interval $f(T_1)$ is much larger than the part of
the interval $T_1'$ which is on the other side of the critical value.
Therefore, the pullback $f^{-n(x)}(D_\phi(T_1))$ is contained in the
lens $D_\phi(T_1)$.
\vglue6pt
Another subcase is the following: suppose that $f^k(x)\in T_1$ ($x\in
(P\cup T_1)\setminus T_2$) and let $k$ be a minimal positive integer
satisfying this property. Put $I(x)=f^{-k}(T_1)$ and $n(x)=k$. Due to
Lemma~\ref{exten_centr} the range of the map $f^k|_{I(x)}$ can be
extended to $T_0$. The pullback of $T_0$ by $f^{-k}$ along the orbit
of $x$ which we denote by $W$, is contained in $\rc(x)$.  Indeed,
suppose that $W\cap S$ is nonempty, so that there is a point $y\in W\cap
S$, and consider two cases.  If $x\in T_1$, then $y\in \partial T_1$
and we would have $f^k(y)\in T_0$ which contradicts the fact that
iterates of the boundary points of $T_1$ never return to the interior
of $T_0$.  On the other hand, if $x\in P$, then $k>m$ because
otherwise we would have $x\notin P$. Now, $f^m(y)$ is either a periodic
point belonging to the boundary of $\Bg_0$ or a point of the forward
orbit of the boundary of $T_1$; thus in any case the point $f^k(y)$
cannot be inside of $T_0$. In both cases we have obtained
contradictions, therefore $W\subset \rc(x)$.
\vglue6pt
By the choice of $T_0$ we know that $W$ contains a $\tau_1$-scaled
neighborhood of $I(x)$, the intervals in the orbit of
$\{f^i(I(x)),\,i=0,\dots k-1\}$ are small and since $I(x)$ is a domain
of the first entry map to $T_1$ the orbit is disjoint.  Hence we can
see that $f^{-k}(D_\phi(T_1))\subset D_{\pi/4}(I(x))\subset
D_\phi(\rc(x))$ (see the choice of the constant $\tau_1$ in the
beginning of the proof).
\vglue6pt
The last case   to consider is the case when $f^i(x)\notin T_1$ for
all $i>0$.  Then $f^i(x)\in \bar P$ for all $i>0$.  Indeed, if
$f^i(x)\not\in \bar P$ for some $i$, then either $f^i(x)\in
[f(\partial T_1),\partial_+ X]$ or $f^{i+j}(x)\in \bar \Bg_0$ for some
$j\leq m$. In the former case we would have $f^{i-1}(x)\in T_1$
(contradiction) and the latter case is impossible because any point of
$\Sigma$ avoids $\Bg_0$. Thus, $x$ belongs to the hyperbolic set
described above, and the sizes of intervals $f^{-i}(\rc(f^i(x)))$ go
to zero as $i\to \infty$. Take $k$ to be so large that $\rc(x)$ is a
$\tau_1$-scaled neighborhood of $f^{-k}(\rc(f^k(x)))$ and
$$\frac{\left|f^{-k}(\rc(f^k(x)))\right|}{|\rc(x)|}<\epsilon.$$
Put $n(x)=k$ and $I(x)=f^{-k}(\rc(f^k(x)))$. By the choice of $m$ we
know that $|\rc(f^k(x))|<C_5$, hence $|f^i(I(x))|<C_3$ for
$i=0,\dots,k$ and\break $\sum_{i=0}^k |f^i(I(x))|<|X|$. As in the previous
case we have $f^{-k}(D_\phi(\rc(f^k(x))))\subset
D_{\pi/4}(I(x))\subset D_\phi(\rc(c))$.

So, we have assigned to each point of $\Sigma'$ one or two
intervals. Now we will show that there are finitely many intervals
of this form whose closures cover all points in $\Sigma'$. First we will
slightly modify these intervals.

When $x\in \Sigma'$, we have assigned to it just one interval which
contains $x$ in its interior. Then we let $\sto{I}(x)$ be the interior of
$I(x)$. Another case: we have assigned to $x$ one interval, say,
$I_-(x)$, but $x$ is its boundary point. Then on the other side of $x$
there is a point $y$ such that the interval $(x,y)$ does not contain
points from the set $\Sigma'$. In this case $\sto I(x)$ is a union of
the interior of $I_-(x)$ and the half interval $[x,y)$. The last case:
there are two intervals assigned to $x$. Let $\sto I(x)$   be the
interior of $I_-(x)\cup I_+(x)$. 

We have covered all points in $\Sigma'$ by open intervals. The set
$\Sigma'$ is compact, therefore there exist finitely many such
intervals which cover $\Sigma'$. Let us denote these
intervals by $\sto I(x_1),\sto I(x_2),\ldots,\sto I(x_N)$. Now,
instead of these intervals consider all the
intervals which are assigned to the points $x_1,\ldots,x_N$,
\ie intervals of the form $I_p(x_i)$, where $p$ is either void or $-$
or $+$ and $i=1,\ldots,N$.
Obviously, the closures of these closed intervals also
cover $\Sigma'$. Moreover, it is easy to see that if the interiors of
two intervals from this set intersect, then one of them is contained
in the other. This is a consequence of the fact that the set $S$
is forward invariant and the boundary points of $I(x)$ are eventually
mapped into $S$. Thus, there exists a finite collection of intervals
of the form $I(x)$ ($I_\pm(x)$) such that the closures of these
intervals cover the whole set $\Sigma'$ and these intervals can
intersect each other only in the boundary points. Denote this
intervals by $I_1,\ldots,I_k$. 

By the construction for each interval $I_i$ there is an integer $n_i$
associated to it. Let $B_i=f^{-n_i}(D_\phi(\rc(f^{n_i}(I_i))))$. We
have the following properties of $I_i$, $n_i$ and $B_i$:
\begin{itemize}
\item $f^{n_i}(I_i)\in \rc$ and $f^{n_i}|_{I_i}$ is monotone if
  $I_i\neq T_2$;
\item If $I_i=T_2$, then $f^{n_i}|_{I_i}$ is unimodal;
\item If $I_i\subset J\in \rc$, 
  then $B_i\subset D_\phi(J)$;
\item If $I_i\neq T_2$, then $B_i\subset
  D_{\pi/4}(I_i)$, thus the domains $B_i$ are disjoint.
\end{itemize}

Let $B=\cup_{i=1}^k B_i$.
It follows that $B$ is a subset of $A$. If $x\in B_i$, put
$F(x)=f^{n_i}(x)$. 

By the very construction of $F$ one can see that it satisfies all the
required properties.\pagebreak \hfill\qed

\section{$C^\omega$ structural stability}\label{ch:Ck}

Here we will prove the $C^k$
structural stability conjecture.

\specialnumber{A} \proclaim{Theorem} \label{AA}
  Axiom {\rm A} maps are dense in the space of $C^\omega(\Delta)$ unimodal
  maps in the $C^\omega(\Delta)$ topology {\rm (}$\Delta$ is an arbitrary
  positive number\/{\rm ).}\/
\endproclaim

We define $C^\omega(\Delta)$ to be the space of real analytic functions
defined on the interval which can be holomorphically extended to a
$\Delta$-neighborhood of this interval in the complex plane.

Let us recall that the map $f$ is {\it regular} if either the
$\omega$-limit set of the critical point does not contain neutral periodic
points or the $\omega$-limit set of $c$ coincides with the orbit of some
neutral periodic point. Any map having negative Schwarzian derivative
is regular. In Section~4.5 we will see that any analytic
map $f$ without neutral periodic points can be included in the family
of regular analytic maps.

\nonumproclaim{Theorem C}
  Let $f_\la: X \hookleftarrow$ be an analytic family of analytic
  unimodal regular maps with a nondegenerate critical point{\rm ,} $\la \in
  \Omega \subset \Rp^N$ where $\Omega$ is an open set. If the family
  $f_\lambda$ is nontrivial in the sense that there exist two maps in
  this family which are not combinatorially equivalent{\rm ,} then Axiom {\rm A}
  maps are dense in this family.  Moreover{\rm ,} let $\Upsilon_{\la_0}$ be
  a subset of $\Omega$ such that the maps $f_{\la_0}$ and $f_{\la'}$
  are combinatorially equivalent for $\la' \in \Upsilon_{\la_0}$ and
  the iterates of the critical point of $f_{\lambda_0}$ do not
  converge to some periodic attractor. Then the set $\Upsilon_{\la_0}$
  is an analytic variety. If $N=1${\rm ,} then $\Upsilon_{\lambda_0}\cap Y${\rm ,}
  where the closure of the interval $Y$ is contained in $\Omega${\rm ,} has
  finitely many connected components.
\endproclaim

{\it Remark}. In Section~4.1 it will be shown that the regularity
condition is superficial if one is concerned only about infinitely
renormalizable maps (or more generally, maps whose $\omega$-limit set
of the critical point is minimal). Thus, the following statements
holds: Let $f_\la: X \hookleftarrow$ be an analytic nontrivial family
of analytic unimodal maps with a nondegenerate critical point, $\la
\in \Omega\break \subset \Rp$, where $\Omega$ is an open set. If the
$\omega$-limit set of the critical point of the map $f_{\lambda_0}$
is minimal, then the set $\Upsilon_{\lambda_0}\cap Y$, where the
closure of the interval $Y$ is contained in $\Omega$, consists of
finitely many points.

In order to underline the main idea of the proof of this theorem we
split it into three parts. First we assume that the map $f$ is infinitely
renormalizable.
In this case the induced quadratic-like map is simpler to study than
the induced polynomial-like map in the other case.  After proving the
theorem in this case we will explain why some extra difficulties in
the general case emerge and then we will show how to overcome them.
Finely we consider the case of Misiurewicz maps (which is the simplest
case).

For the reader's convenience we collect all theorems about
quasi-conformal maps which we will use intensively  in
Appendix~\ref{sec:appendix}.

\demo{{\rm 4.1.} The case of an infinitely renormalizable map} 
In this section we will proof the following lemma: \enddemo
\specialnumber{4.1}\proclaim{Lemma}\label{lm:C-minimal}
  Let $f_\la: X \hookleftarrow$ be an analytic family of analytic unimodal maps
  with a nondegenerate critical point{\rm ,} $\la \in \Omega \subset \Rp^N$ where $\Omega$ is
  a open set. Suppose that the map $f_{\lambda_0}$ is infinitely
  renormalizable. Then there is a neighborhood $\Omega'$ of $\lambda_0$ such that
  the set $\Upsilon_{\lambda_0}\cap \Omega'$ is an analytic variety.
\endproclaim

This lemma remains true if instead of assuming that the map $f_{\lambda_0}$ is
infinitely renormalizable, we assume that the $\omega$-limit set of the
critical point of this map is minimal.  Note that we do not assume
here that the family $f$ is regular.

We can assume that $\lambda_0=0$.
 
From Theorem~\ref{thr:ren_minimal} we know that if the map is analytic
and infinitely renormalizable, then there is an induced quadratic-like
map $F_0:B \to A$, where $B\subset A \subset \Cp$ are simply connected
domains and the modulus of the annulus $A\sm B$ is not zero.

The map $F_0$ is the extension of some iterate of the map $f_0$ to the
domain $B$, i.e.,  $F_0|_B=f_0^n$.
If we take a small neighborhood $D\subset\Cp^N$ of $0$ in the parameter
space, then the map $F_\la=f_\la^n$ will have the extension to some
domain which contains $B$ for any $\la\in D$. Fix the domain $A$ and let
$B_\la$ be a preimage of the domain $A$ under the map $F_\la$ where
$\la \in D$ and let $B_\la \subset A$.

Define the map $\phi_\la:\partial B_0 \cup \partial A \to \partial B_\la \cup \partial A$ by the following
formula: $\phi_\la(z)=F_\la^{-1}\circ F_0(z)$ where $\la \in D$, $z \in \partial B_0$
and $\phi_\la(z)=z$ for $z\in \partial A$. The map $F_\la$ is not invertible, but
if $\phi$ is continuous with respect to $\la$ and $\phi_0=\mbox{id}$, then
it is defined uniquely.

For fixed $z$ the map $\phi_\la(z)$ is holomorphic with respect to $\la$.
Shrinking the neighborhood $D$ if necessary, we can suppose that the
map $z \mapsto \phi_\la(z)$ is injective for fixed $\la\in D$.  Due to $\lambda$-lemma
(Theorem~\ref{la-lemma}) the map $\phi_\lambda$ can be extended to the annulus
$A\setminus B_0$ in the q.c. (quasiconformal) way. Denote this extension by $h_\lambda^0:A\setminus
B_0\to A\setminus B_\lambda$.  Thus, $h_\lambda^0$ is a q.c.~homeomorphism and its Beltrami
coefficient $\nu_\lambda^0$ is a holomorphic function with respect to $\lambda\in D$.

Denote the pullback of the Beltrami coefficient $\nu_\la^0$ by the map
$F_0$ as $\nu_\la$; i.e.,  if $F_0^k(z)\in A\sm B$, then
$\nu_\la(z)=F_0^{k\,*} \nu_\la^0(F_0^k(z))$. On the filled Julia set
of $F_0$ and outside of the domain $A$ we set $\nu_\la$   equal to~0. It is easy to see that since $\la
\mapsto \nu_\la^0(z)$ is analytic the map $\la \mapsto \nu_\la(z)$ is analytic as well.

According to the measurable Riemann mapping Theorem~\ref{mRm} below, there is
a family   q.c.~homeomorphism  $h_\la:\Cp\to \Cp$ whose Beltrami
coefficient is $\nu_\la$ and which is normalized such that
$h_\la(\infty)=\infty$, $h_\la(a^-)=a^-$, $h_\la(a^+)=a^+$ where the
$a^\pm$ are two points of the intersection of $\partial A$ and the real
line.

Since the map $F_0$ conserves the Beltrami coefficient $\nu_\la$ the
map
$$G_\la=h_\la \circ F_0 \circ h_\la^{-1}: B_\la \to A$$
is
holomorphic. Due to the Ahlfors-Bers Theorem~\ref{ahl_be} the map $\la
\mapsto G_\la(z)$ is analytic for the fixed point $z$. Thus $G$ is an
analytic family of holomorphic quadratic-like maps.

\specialnumber{4.2}\proclaim{Lemma} \label{la-1}
  The maps $f_0$ and $f_\la$ are combinatorially equivalent if and
  only if $F_\la=G_\la$.
\endproclaim

\begin{plm} 
If $F_\lambda=G_\lambda$, then $F_\lambda$ and $F_0$ are topologically
conjugate; hence $f_\lambda$ and $f_0$ are combinatorially
equivalent. 

If $f_0$ and $f_\la$ are combinatorially equivalent, then the
maps $F_0$ and $F_\la$ are combinatorially equivalent as well. Due to
the rigidity theorem and straightening Theorem~\ref{straigh} we
know that there is a q.c.~homeomorphism $\tilde H:\Cp\to\Cp$ which is
a conjugacy between $F_0$ and $F_\la$ on their Julia sets; i.e.,
$\tilde H\circ F_0|_J=F_\la \circ \tilde H|_J$ where $J$ is the Julia
set of the map $F_0$.

Define a new q.c.~homeomorphism $H^0$ in the following way:
$$
H^0(z)=\left\{
\begin{array}{ll}
z & \mbox{if } z\notin A\\
h^0_\la(z) & \mbox{if } z \in A\sm B\\
\tilde H(z) & \mbox{if } z \in B(J)
\end{array}
\right.
$$
where $B(J)$ is a neighborhood of the Julia set $J$ such that $B(J)
\subset B$. In the annulus $B\sm B(J)$ the q.c.~homeomorphism $H^0$ is
defined in an arbitrary way.

Consider the sequence of q.c.~homeomorphisms $H^i$ which are defined
by the formula $H^{i+1}=F_\la^{-1} \circ H^i \circ F_0$. The map
$F_\la$ is not invertible, but $H^{i+1}$ is defined correctly because of
the homeomorphism $\tilde H$ and as a consequence the homeomorphism
$H^i$ maps the orbit of the critical point of $F_0$ onto the orbit of
the critical point of $F_\la$. Since the maps $F_0$ and $F_\la$ are
holomorphic the distortion of $H^{i}$ does not increase with $i$. So
the sequence $\{H^i\}$ is normal and we can take a subsequence
convergent to some limit $\hat H$ which is also a q.c.~homeomorphism.
Taking a limit in the equality $H^{i+1}=F_\la^{-1} \circ H^i \circ
F_0$ we obtain that the homeomorphism $\hat H$ is a conjugacy between
$F_0$ and $F_\la$; i.e.,  $F_\la \circ \hat H=\hat H \circ
F_0$. On the other hand, it is easy to see that the Beltrami
coefficient of $\hat H$ coincides with the Beltrami coefficient
$\nu_\la$. Indeed, outside of $A$ both coefficients are zero. In the
domain $A\sm J$ both coefficients are obtained by pulling back the
Beltrami coefficient $\nu_\la^0$. On the Julia set the Beltrami
coefficient of $\hat H$ is equal to the Beltrami coefficient of
$\tilde H$ which is 0 because of the rigidity theorem. The
homeomorphism $\hat H$ is normalized in the same way as $h_\la$, so that by
the measurable Riemann mapping theorem these homeomorphisms coincide.
From the very definition of the map $G_\la$ we obtain that
$F_\la=G_\la$.  \end{plm}


Due to the previous lemma $f_0$ and $f_\lambda$ are combinatorially
equivalent if and only if $F_\lambda=G_\lambda$. So, the solution with respect to
$\lambda$ of the equation $F_\lambda=G_\lambda$ is the set $\Upsilon_0\cap D$. Since this equation
is holomorphic, its solution is an analytic variety.

\demo{{\rm 4.2.} The case of a finitely renormalizable {\rm (}\/nonrenormalizable\/{\rm )} map}
In the previous section the domain $A\sm B$ had the nice boundary
which was a union of two Jordan curves. In the general case this is
false. Indeed, recall the structure of the domains $A$ and $B$ which
is given in Section~\ref{ch:p_like}. The domain $A$ is a union of
finitely many lenses based on the real line. Inside of each lens there
are finitely many quasilenses which are connected components of the
domain $B$ (see Figure~2). Thus, if $ A^{x_0}\subset A$ is a
connected component of the domain $A$, then the set $A^{x_0}\sm B$
consists of 1 or 2 connected components which can have cusps or angles
on their boundaries (recall that $A^x$ denotes a connected component
of A containing the point $x$).

Notice that the family $f_\lambda$ consists of regular maps so that we will
not have neutral periodic points on the boundary of the domains $A$
and $B$.

Let $a$ be a periodic point from the set
$E=\omega(\partial(A\cap\Rp))$ (see   \S 3.2).
For   simplicity we will assume that the point $a$ is a fixed point.
Denote the multiplier of the map $F_\la$ at the point $a$ as $d_\la$
and let $\partial A^{x_0}$ and $\partial B^{x_0}$ contain the point~$a$. If on the boundary of the domain $
A^{x_0}$ we define the map
$h_\la^0$ to be the identity, then on the boundary of the domain $B$
near the point $a$ we will have $h_\la^0(z)=d_0/d_\la \,z +\cdots$
because the map $h_\la^0$ has to conjugate the maps $F_0$ and
$F_\lambda$ on the boundary of $B$; i.e.,   $h_\la^0|_{\partial
  A}\circ F_0|_{\partial B}=F_\la|_{\partial B_\la} \circ
h_\la^0|_{\partial B}$. At the point $a$ the boundaries of the domains
$B$ and $A$ are tangent to each other, and if the multiplier $d_\la$
changes with $\la$, then the derivative of $h_\la^0$ in the direction
of $\partial A$ is 1 and in the direction of $\partial B$ is
$d_0/d_\la$. One can easily check that a homeomorphism $h_\la^0$
defined on the domain $A\sm B$ cannot be quasiconformal.

As a result of this discussion we conclude that we have to deform the
domain $A_\la$ as well in order to construct the q.c. homeomorphism
$h_\la^0$.

Now we will prove Lemma~\ref{lm:C-minimal} in the case when the map
$f_0$ is finitely  renormalizable. 

\specialnumber{4.3}\proclaim{Lemma}\label{lm:C-nonminimal}
  Let $f_\la: X \hookleftarrow$ be an analytic regular family of analytic unimodal
  maps with a nondegenerate critical point{\rm ,} $\la \in \Omega \subset \Rp^N$ where
  $\Omega$ is an open set. Suppose that the map $f_{\lambda_0}$ is finitely
  renormalizable. Then there is a neighborhood $\Omega'$ of $\lambda_0$ such that
  the set $\Upsilon_{\lambda_0}\cap \Omega'$ is an analytic variety.
\endproclaim

Recall the notation   used in Section~3.2. According to
Theorem~\ref{p_like}, for our map $f_0$ there is an induced
polynomial-like map $F_0:B_0\to A_0$.  The set $S$ consists of points
where the domain $A_0$ has singularities. This set is finite and
forward invariant, so that it has periodic points and let $E$ denote this
subset of periodic points. Any point from the set $S$ is mapped into
$E$ after some iterations.

We can make an analytic change of the coordinate which also depends on
the parameter $\lambda$ analytically in such a way that the set $S$
does not move with the parameter $\lambda$ for small $\lambda$. So in
this section we will assume that the set $S$ does not depend on
$\lambda$.

Take any periodic point $r$ from the set $E$ and let $m$ be the
period of this periodic point $r$. Let $x$ be a local coordinate in
the neighborhood of the point $r$ and let the map $f_\lambda^{m}$
have the following series expansion:
$$
  f_\lambda^m(x)   =  d_\lambda x+q_\lambda x^2+\Oi(x^3).
$$
The coefficients $d_\lambda$ and $q_\lambda$ depend analytically on the
parameter $\lambda$.

Our goal is the construction of a q.c.~homeomorphism
$h_\lambda^0:A_0\sm B_0\to \Cp$ which conjugates the maps $F_0$ and
$F_\lambda$ on the domain $A_0\setminus B_0$.

Assume that $d_\lambda>0$ and let $ A^{x_0}_0\supset  B^{x_0}_0$ be connected
components of the domains $A_0$ and $B_0$ which have the point $r$ in
their boundaries. It follows from the construction of the domains $A_0$
and $B_0$ that at the point $r$ the boundaries of $A^{x_0}$ and $B^{x_0}$
are tangent to each other and that this tangency is quadratic.  We
will look for the map $h_\lambda^0$ near the point $r$ in the
following form:
$$
h_\lambda^0(z)=(z-r)^{l_\lambda} b_\lambda(z-r)(1+\oi(z-r)),
$$
where $b(z)$ is a holomorphic function such that $b(0)\neq 0$.

Since the map $h_\lambda^0$ should conjugate the maps $F_0$ and
$F_\lambda$ we obtain the following equation for $h_\lambda^0$:
$$
h_\lambda^0\circ f_0^m=f_\lambda^m \circ h_\lambda^0.
$$

Solving this equation we obtain the series expansion of $h_\lambda^0$:
$$
h_\lambda^0(z)=(z-r)^{l_\lambda}+ \alpha_\lambda (z-r)^{2l_\lambda}+
\beta_\lambda (z-r)^{l_\lambda+1}+\Oi((z-r)^\kappa)
$$
where
\begin{eqnarray*}
  l_\lambda&\hskip-6pt =\hskip-6pt&\frac{\ln(d_\lambda)}{\ln(d_0)}, \qquad 
  \alpha_\lambda = \frac{q_\lambda}{d_\lambda^2}\\
  \beta_\lambda&\hskip-6pt=\hskip-6pt&\frac{l_\lambda q_0}{d_0(1-d_0)},\qquad
  \kappa =  \min(3l_\lambda,2l_\lambda+1).
\end{eqnarray*}

Now to each point of the set $S$ we associate a jet by the following
rule: first, from each periodic orbit of the set $E$ take a
representative and denote this set of representatives as $E'$. For a
point $r\in E'$ the corresponding jet $j_{r,\lambda}$ is defined as
$x^{l_\lambda}+ \alpha_\lambda x^{2l_\lambda}+ \beta_\lambda
x^{l_\lambda+1}+\Oi(x^\kappa)$ where $l_\lambda$, $\alpha_\lambda$ and
$\beta_\lambda$ are calculated according to the formulas above. If
$a\in S\setminus E'$, then some iteration of $a$ is mapped into the
set $E'$, so that $f^n(a)=r$ where $r$ is some element of the set~$E'$.
Then at the point $a$ the jet $j_{a,\lambda}$ is defined as
$f_\lambda^{-n}\circ j_{r,\lambda} \circ f_0^n$. Certainly, we
truncate the terms of order $\Oi(x^{\kappa})$ and higher.

Now, at each point of the set $S$ we have a jet which depends on the
parameter $\lambda$.

The family of maps $\phi_\lambda:\partial A_0\cup \partial B_0 \to
\Cp$ will be defined first on the boundary of the domain $A_0$.  Let it
satisfy the following conditions:
\begin{itemize}
\item $\phi_0=\mbox{id}$;
\item For fixed $z\in \partial A$ the map $\lambda\mapsto
  \phi_\lambda(z)$ is analytic;
\item For fixed $\lambda$ the map $z\mapsto \phi_\lambda(z)$ is
  differentiable and nonneutral for $z\in \partial A_0\setminus
  S$;
\item For any $r\in S$ we have
  $\phi_\lambda(z)=j_{r,\lambda}(z-r)+\Oi((z-r)^\kappa)$. 
\end{itemize}
One can easily construct the map $\phi_\lambda$ satisfying these
conditions. 
\enddemo
On the boundary of the domain $B_0$ we define the map $\phi_\lambda$
in such a way that $\phi_\lambda$ conjugates the maps $F_0$ and
$F_\lambda$; i.e.,
$$
\phi_\lambda|_{\partial A} \circ F_0|_{\partial B}
=F_\lambda|_{\partial B_\lambda} \circ \phi_\lambda|_{\partial B}.
$$
Thus
$$
\phi_\lambda|_{\partial B_0}=F_\lambda^{-1}|_{\partial A_\lambda} \circ
\phi_\lambda|_{\partial A_0} \circ F_0|_{\partial B_0}
$$
where $\partial A_\lambda=\phi_\lambda(\partial A_0)$.

From the construction it follows that at the points where the domain
$A_0\setminus B_0$ has quadratic singularities (i.e.\  at
points of the set $S$) we have
$$
\phi_\lambda(z-a)=\gamma_{a,\lambda} (z-a)^{l_\lambda}
+\alpha_{a,\lambda}(z-a)^{2l_\lambda}
+\beta_{a,\lambda}(z-a)^{l_\lambda+1}
+\Oi((z-a)^\kappa)
$$
where $a\in S$ and $z\in \partial A_0 \cup \partial B_0$.

\figin{fig4}{800}

\begin{quote} Figure 4. A connected component of the domain
  $A_0$. At the point $b$ the angle is not zero.
\end{quote}

If $b$ is a singularity of the domain $A \setminus B$ where
this domain has a nonzero angle (i.e.\  $b$ is a point of the intersection
of the closure of two connected components of the domain $B_0$), denote two
arcs which are boundary arcs of the domain $B$ and which intersect at
$b$, as $J^-$ and $J^+$ (see Fig.\ 4). Let $F_0|_{J^i}=f^{k_i}$
for $i=-,+$. The numbers $k_-$ and $k_+$ do not necessarily coincide.
Therefore, the jets of the maps $\phi_\lambda|_{J^-}$ and
$\phi_\lambda|_{J^+}$ are different. However, the exponents of the
leading terms of these jets do coincide. So, in the neighborhood of
the point $b$ we have
$$
\phi_\lambda(z)=\gamma_{i,\lambda}
(z-b)^{l_\lambda}(1+\Oi((z-b)^{\min(l_\lambda,1)})) 
$$
for $z\in J^i$ where $i=-,+$, and $\gamma_{i,\lambda}$ is holomorphic
with respect to $\lambda$, $\gamma_{i,0}\neq 0$ and
$\gamma_{i,\lambda}$ is real for real $\lambda$.

\specialnumber{4.4}\proclaim{Lemma}
  There is a small neighborhood $D\subset \Cp^N$ of 0 such that for fixed
  $\lambda\in D$ the map $\phi_\lambda:A_0\setminus B_0 \to \Cp$ defined above is injective.
\endproclaim
 
\begin{plm}
First, we will check that the map $\phi_\lambda$ is injective in some
small neighborhood of the point $b$ where we have a nonzero angle. 

Let $x$ be a local coordinate in the neighborhood of $b$ and let the
curves $J^-$ and $J^+$ have the parametrizations $x=u_- t +\Oi(t^2)$
and $x=u_+ t+\Oi(t^2)$, where $t\in \Rp$ and $u_-,u_+\in \Cp$. Since
the angle at $b$ is nonzero the ratio $\frac{u_-}{u_+}$ cannot be
real. 

Suppose that $\phi_\lambda$ is not injective. Then there are real
numbers $t_-$ and $t_+$ such that 
$$
\gamma_{-,\lambda}\, u_-^{l_\lambda}\, t_-^{l_\lambda} \,
(1+\Oi(t_-^{\min(l_\lambda,1)}))=
\gamma_{+,\lambda}\, u_+^{l_\lambda}\, t_+^{l_\lambda} \,
(1+\Oi(t_+^{\min(l_\lambda,1)})).
$$

For small $\lambda$ the exponent $l_\lambda$ is close to 1. Hence, for
small $\lambda$ the imagery part of
$\frac{\gamma_{-,\lambda}}{\gamma_{+,\lambda}}
\left(\frac{u_-}{u_+}\right)^{l_\lambda}$ is bounded away from
0. Thus, for small $\lambda$ and $t_-$, $t_+$ the equation
$$
\frac{\gamma_{-,\lambda}}{\gamma_{+,\lambda}}
\left(\frac{u_-}{u_+}\right)^{l_\lambda}=
\left(\frac{t_-}{t_+}\right)^{l_\lambda}(1+\Oi(t_-^{\min(l_\lambda,1)})+
\Oi(t_+^{\min(l_\lambda,1)}))
$$
does not have real solutions.

Consider now the point $a\in S$ where we have a quadratic singularity.
Let us again parametrize the boundaries of $A_0$ and $B_0$ in the
neighborhood of $a$ by $x=ut+v_-t^2+\Oi(t^3)$ and
$x=ut+v_+t^2+\Oi(t^3)$ where $u$ is a complex number, $v_-,\,v_+$ are
real numbers and $v_-\neq v_+$.

The equation we have to solve is the following:
\begin{eqnarray*}
&&\hskip-22pt\gamma_\lambda\,(ut_-+v_-t_-^2)^{l_\lambda}+
\alpha_\lambda\,(ut_-+v_+t_-^2)^{2l_\lambda}+
\beta_\lambda\,(ut_-+v_-t_-^2)^{l_\lambda+1} +
\Oi(t_-^{\kappa_\lambda})\\[6pt]
&&  \!\! =
\gamma_\lambda\,(ut_++v_-t_+^2)^{l_\lambda}+
\alpha_\lambda\,(ut_++v_+t_+^2)^{2l_\lambda}+
\beta_\lambda\,(ut_++v_-t_+^2)^{l_\lambda+1} +
\Oi(t_+^{\kappa_\lambda}).
\end{eqnarray*} 

After simplification we obtain:\eject
\begin{eqnarray*}
\noalign{\vskip-18pt}
&&\hskip-22pt\gamma_\lambda(ut_-)^{l_\lambda}+
\gamma_\lambda l_\lambda u^{l_\lambda-1}v_- t_-^{l_\lambda+1}+
\alpha_\lambda(ut_-)^{2l_\lambda}+
\beta_\lambda(ut_-)^{l_\lambda+1} +
\Oi(t_-^{\kappa_\lambda})\\[6pt]
&&\quad= 
\gamma_\lambda(ut_+)^{l_\lambda}+
\gamma_\lambda l_\lambda u^{l_\lambda-1}v_+ t_+^{l_\lambda+1}+
\alpha_\lambda(ut_+)^{2l_\lambda}+
\beta_\lambda(ut_+)^{l_\lambda+1} +
\Oi(t_+^{\kappa_\lambda}).
\end{eqnarray*}

One can easily see that this equality implies that
$t_-=t_+ +\frac{v_+-v_-}{u} t_+^2\break+ \oi(t_+^2)$. However, $v_+-v_-$ is a
real number and $u$ is complex, so if $t_+$ is a small real number,
then $t_-$ is complex. Thus, for small $\lambda$ the map
$\phi_\lambda$ is injective in small neighborhoods of the singular
points.

If at some point  the boundary of $B_0$ or $A_0$ is smooth, then for
small $\lambda$ the map $\phi_\lambda $ is injective as well in some
neighborhood of this point. By compactness arguments we obtain that
for small $\lambda $ the map $\phi_\lambda$ is injective.
\end{plm}

According to the $\lambda$-lemma we can extend the map $\phi_\lambda$
to the domain\break $A_0\setminus B_0$.  In other words, there is a family
of q.c.~homeomorphisms\break $h_\lambda^0:A_0\setminus B_0\to \Cp$ where
$\lambda$ is in some small neighborhood of the point 0. This family
satisfies the following conditions:
\begin{itemize}
\item $h_0^0=\mbox{id}$;
\item For the fixed parameter $\lambda$ the map $h_\lambda^0$ is a
  q.c.~homeomorphism and $h_\lambda^0|_{\partial A \cup \partial
  B}=\phi_\lambda$;
\item For fixed $z\in A_0\setminus B_0$ the maps $\lambda\mapsto
  h_\lambda^0(z)$ and $\lambda\mapsto \nu_\lambda^0(z)$ are analytic
  where $\nu_\lambda^0$ is the Beltrami coefficient of $h_\lambda^0$.
\end{itemize}

Now we have the map $h_\lambda^0$, so we can construct the
q.c.~homeomorphism $h_\lambda$ and the analytic family $G$.
Lemma~\ref{la-1} still holds, but we have to alter its proof because
we cannot use the straightening theorem any more. Instead of it we will
use the following theorem (see \cite{GS}, \cite{Lyu4}).

\proclaim{Theorem}
  Let $R_0:\hat B_0\to \hat A_0$ and $R_1:\hat B_1 \to \hat A_1$ be 
  holomorphic box mappings such that $R_0$ and $R_1$ are
  combinatorially equivalent{\rm ,} and the moduli of the annuli $\hat
  A_i\setminus \hat B_i^x$ are uniformly bounded away from zero for
  all $x\in \hat B_i\cap\Rp${\rm ,} where $ B_i^x$ is a connected component
  of $\hat B_i$ containing the point $x${\rm ,} $i=0,1$.  Moreover{\rm ,} suppose
  that there is a quasisymmetric  homeomorphism $Q$ such that $Q\circ
  R_0|_{\partial \hat B_0 \cap \Rp}=R_1\circ Q|_{\partial \hat B_1
    \cap \Rp}$. Then the maps $R_0$ and $R_1$ are {\rm q.c.}~conjugate on
  their postcritical sets.
\endproclaim 

Consider the map $F_0:B_0\to A_0$ which is induced by the map $f_0$.
Let $A_0^c$ be a connected component of $A_0$ which contains the
critical point. If $ B_0^{x}$ is a connected component of $B_0$ which
is mapped onto $A_0^c$ by $F_0$ (this is equivalent to saying that
$F_0(x)\in A_0^c$), then the domain $ B_0^x$ is disjoint from the
boundary of the domain $A_0$ (see Theorem~\ref{p_like}). Since there
are only finitely many connected components of the domain $B_0$ we see
that there is a positive number $C_6$ such that
$\mathop{\mathrm{mod}}( A_0^x \setminus B_0^x)>C_6$ for any $x\in
B_0\cap \Rp$ such that $F_0(x)\in A_0^c$.

Denote the first return map of the map $F_0$ to the domain $A_0^c$ by
$R_0$ and $A_0^c$ by $\hat A_0$. It is easy to see that $R_0$ is a
holomorphic box mapping and that the moduli of the annuli $\hat
A_0\setminus \hat B_0^x$ with $x\in \hat B_0\cap\Rp$ are uniformly
bounded away from zero by the constant $C_6$, where $\hat B_0^x
\subset \hat A_0$ is a connected component of the domain of definition
of the map $R_0$.


In a similar way we can define the first entry map $R_\lambda$. In
order to apply the previous theorem to the maps $R_0$ and $R_\lambda$
and to find the q.c.~conjugacy between $R_0$ and $R_\lambda$ on their
postcritical sets we have to construct the q.s.~homeomorphism $Q$. It
is easy to do using the following observations: first, the maps $f_0$
and $f_\lambda$ are regular, hence they have no neutral periodic
points (we have supposed that the critical points are recurrent); in
this case the set of points which do not belong to the basin of
attraction and whose iterates do not enter some neighborhood of the
critical point is a hyperbolic set; since the maps $f_0$ and
$f_\lambda$ are conjugate these corresponding hyperbolic sets are
conjugate as well and this conjugacy $Q$ is quasi-symmetric. This can
be proved using the same ideas as for the Misiurewicz maps; see, for
example, \cite{dMvS}.  Obviously, the set $\partial \hat B_0\cap\Rp$
is a subset of the hyperbolic set which consists of points whose
iterates do not enter the interval $\hat B^c_0\cap\Rp$ (and do not
belong to the basin of attraction).  Another way to see the existence
of this q.s.~homeomorphism in our case is the following: the set
$\partial\hat B_\lambda\cap \Rp$ consists of preimages of points in
$\partial B_\lambda\cap \Rp$ and it varies holomorphically with respect
to $\lambda$. Moreover, this set is a part of some hyperbolic set,
hence it persists for small $|\lambda|$ (even if $\lambda$ is
complex). Now we can apply the  $\lambda$-lemma and get a q.c.~homeomorphism
which maps $\partial\hat B_0\cap \Rp$ onto $\partial\hat B_\lambda\cap \Rp$.

According to the previous theorem there is a q.c.~homeomorphism $H$
which conjugates the maps $R_0$ and $R_\lambda$ on their postcritical
sets if the maps $f_0$ and $f_\lambda$ are conjugate. By pulling
forward we can find a q.c.~homeomorphism $\tilde H$ which is a
conjugacy of the maps $F_0$ and $F_\lambda$ on their postcritical
sets.

Having this map $\tilde H$ we can proceed with the proof exactly in
the same way as in Section~4.1. Indeed, we
can construct a sequence of q.c.~homeomorphisms $H^k$ and take a
subsequence converging to $\hat H$. If the map $F_0$ is
nonrenormalizable, then the Julia set of $F_0$ has zero Lebesgue
measure. The proof of this fact is given in
Appendix~5.4. Thus, we can again conclude that
$h_\lambda=\hat H$ and therefore $G_\lambda=F_\lambda$ if $F_0$ is combinatorially
equivalent to $F_\lambda$.

\vglue6pt
 4.3. {\it The case of a Misiurewicz map}.
Finally, let us consider the case when $f_0$ is a Misiurewicz map.

\specialnumber{4.5}\proclaim{Lemma}\label{lm:C-misiurewicz}
  Let $f_\la: X \hookleftarrow$ be an analytic regular family of analytic unimodal
  maps with a nondegenerate critical point{\rm ,} $\la \in \Omega \subset \Rp^N$ where
  $\Omega$ is an open set. Suppose that the map $f_{\lambda_0}$ does not satisfy
  Axiom {\rm A} and that the critical point of $f_{\lambda_0}$ is nonrecurrent.
  Then there is a neighborhood $\Omega'$ of $\lambda_0$ such that the set
  $\Upsilon_{\lambda_0}\cap \Omega'$ is an analytic variety.
\endproclaim

Since the critical point of $f_0$ is nonrecurrent,   there exists a
neighborhood $U$ of $c_0$ such that $f_0^n(c_0)\not \in U$ for all
$n>0$, where $c_0$ is a critical point of $f_0$.  Let $\Sigma_0$ be a set
of points which do not belong to the basin of attraction and whose
forward orbits under iterates of $f_0$ do not enter $U$. Obviously,
$f_0(c_0)\in \Sigma_0$ which is a closed set and does not contain neutral
periodic points because $f_0$ is a regular map and we have assumed
that the iterates of the critical point do not converge to a periodic
attractor. Due to Ma\~{n}\`e's theorem $\Sigma_0$ is a hyperbolic set and there
exists a neighborhood $D\subset \Cp^N$ of $0$ such that when $\lambda$ is in $D$,
$f_\lambda$ has a hyperbolic set $\Sigma_\lambda$ close to $\Sigma_0$ and the dynamics of
$f_\lambda$ on $\Sigma_\lambda$ is conjugate to the dynamics of $f_0$ on $\Sigma_0$.  Thus
there exists a homeomorphism $h_\lambda:\Sigma_0\to \Sigma_\lambda$.  The set $\Sigma_\lambda$ depends
holomorphically on $\lambda$. Indeed, the periodic points in $\Sigma_\lambda$ depend
holomorphically on $\lambda$ and they are dense in $\Sigma_\lambda$.  Applying the
$\lambda$-lemma we can conclude that for fixed $z$ the map $h_\lambda(z)$ is
holomorphic.

The maps $f_0$ and $f_\lambda$ are combinatorially equivalent for some
$\lambda\in D\cap \Rp^N$, if and only if
$h_\lambda(f_0(c_0))=f_\lambda(c_\lambda)$. The last equation is
analytic with respect to $\lambda$; hence its solution is an analytic
variety.

\vglue6pt
  4.4. {\it Density of Axiom {\rm A} in regular families}.
Now we finish the proof of Theorem~\ref{real_anal}.

First let us consider the case $N=1$.  Suppose that $f_{\lambda_0}$ does not
satisfy Axiom A and that the set $\Upsilon_{\lambda_0}$ contains infinitely many
points in $Y$. Since $\Upsilon_{\lambda_0}$ is an analytic variety,   it is an open
set. However, from kneading theory we know that this set of
combinatorially equivalent maps should be closed. We have arrived at
a contradiction and hence the set $\Upsilon_0$ has only finitely many
points.
\pagegoal=50pc

Now we shall prove that Axiom A maps are dense in $\Omega$. We have already
shown that if the iterates of the critical point of some map $f_{\lambda_0}$
do not converge to a periodic attractor, then one can perturb this map
within the family $f_\lambda$ to some other map which is not combinatorially
equivalent to $f_{\lambda_0}$. The kneading invariant changes continuously
with $\lambda$; hence there is a map $f_{\lambda_1}$ in the family close to
$f_{\lambda_0}$ such that the iterates of its critical point converge to
some periodic attractor. If this attractor is hyperbolic, we are done
because then there are no neutral periodic orbits and the map is an
Axiom A map. The other case is that the attractor is a neutral
periodic orbit. The multiplier of this periodic orbit is an analytic
function with respect to $\lambda$;  hence either there are maps in the
family $f_\lambda$ close to $f_{\lambda_1}$ which do not have a neutral periodic
orbit of the same period or such a neutral periodic orbit exists for
all $\lambda\in \Omega$. In the former case we can find a map close to $f_{\lambda_1}$
such that the iterates of its critical point converge to a hyperbolic
periodic orbit (this orbit appears after a bifurcation of the neutral
periodic orbit), and this map is an Axiom A map. The latter case is
impossible because in this case the iterates of the critical point
should converge to this neutral periodic orbit for all maps in the
family and hence all maps in the family would be combinatorially
equivalent.
 \pagebreak

4.5. {\it Construction of a regular family}.
Now we are going to show how to derive Theorem~\ref{AA} from
Theorem~\ref{real_anal} and first we will study some properties of
regular maps.
 \pagegoal=48pc
\specialnumber{4.6}\proclaim{Lemma}
   Any regular map $f\in C^3$ with a recurrent critical point
  has its neighborhood in the space of $C^3$ unimodal maps
  consisting of regular maps.
\endproclaim

\begin{plm} Since $f$ is regular and its critical point is recurrent, the map
$f$ has no neutral periodic points. Consider a nice interval $I_f$
around the critical point such that the first return map to $f(I_f)$
has negative Schwarzian derivative (see Theorem~\ref{nsch}). It can be
easily shown that if a map $g$ is $C^3$ close to $f$, then for the map
$g$ there is a nice interval $I_g$ close to $I_f$ such that the first
return map of $g$ to $g(I_g)$ has negative Schwarzian derivative as
well. Let $J$ be an interval containing the critical point and let the
interval $I_f$ strictly contain $J$. The set of points whose iterates
under the map $f$ never enter the interval $J$ is a union of some
hyperbolic set, periodic attractors and points whose iterates converge
to the periodic attractors. If $g$ is $C^3$ close enough to $f$, then
the interval $I_g$ will contain $J$ and the hyperbolic set and its periodic attractors persist. In this case
the map
$g$ is regular. Indeed, if $g$ has a neutral periodic point, then the orbit of this
point necessarily passes through the interval $g(I_g)$. The first
return map of $g$ to $g(I_g)$ has negative Schwarzian derivative, and,
hence, iterates of the critical point have to converge to this neutral
point (this is a standard fact, see \cite{Sin}).  \end{plm}


The set of unimodal maps in $C^\omega(\Delta)$ which have a neutral
periodic orbit of period $K$ is an analytic variety of codimension $1$;
thus the complement of this set is open and dense in
$C^\omega(\Delta)$. The set of maps which do not have neutral periodic
orbits is equal to the intersection of all such   complements for
$K=1,2,\dots$. Due to the Baire theorem this set is dense in
$C^\omega(\Delta)$ as well. Thus we have proved the following lemma:

\specialnumber{4.7}\proclaim{Lemma} \label{lm:regular}
  The set of regular maps is dense in the space of unimodal maps
  $C^\omega(\Delta)$.
\endproclaim

{\it Proof of Theorem} A.
We will show that any regular map with a recurrent critical
point can be included in a nontrivial analytic family of regular
analytic unimodal maps. This will imply Theorem~\ref{AA}. Indeed,
since the regular maps are dense in $C^\omega(\Delta)$ we can first
perturb the given map to a regular map, and then we can construct a
nontrivial family of regular analytic maps and apply
Theorem~\ref{p_like}.

First   notice that if the map we need to perturb is infinitely
renormalizable, then we can take any nontrivial family passing
through this map and apply the statement formulated in the remark
after Theorem~\ref{real_anal}; see also Section~4.1. In
this way we can obtain a map close to the original map such that the
iterates of its critical point converge to some periodic attractor. If
this map has neutral periodic points, it is easy to perturb it to a
map which does not have neutral periodic orbits, however the iterates
of its critical point still converge to a periodic hyperbolic
attractor.  Obviously, this will be an Axiom A map.  The same
arguments apply to the case when the map we need to perturb is a
Misiurewicz map. Indeed, in Section~4.3 we have only
used the regularity of the map $f_0$ itself and we have never used the
regularity of other maps in the family. Thus, we have only to
construct a perturbation of an analytic unimodal regular
nonrenormalizable map with a quadratic recurrent critical point.



Now we are going to construct a perturbation of $f$. First, it will be
only a $C^3$ perturbation. 

For any $\epsilon>0$, Theorem~3.1 gives   a
polynomial-like map $F_\epsilon:B_\epsilon\to A_\epsilon$ induced by
$f$.  Let $A_\epsilon^{c}$ be a connected component of $A_\epsilon$
containing the critical point and let
$a_\epsilon=f(\partial(A_\epsilon^c\cap \Rp))$.  The interval
$f(B_\epsilon^c)\cap \Rp$ has two boundary points as well and we let
$b_\epsilon$ be one of these boundary points which does not have two
real preimages under $f$. Just to fix the notation let us assume that
$a_\epsilon<b_\epsilon$, which corresponds to the case when the map $f$
first increases and then decreases.

Let the function $p_{\epsilon,\lambda}:\Rp\to\Rp$ be given by the following
formula:
$$
p_{\epsilon,\lambda}(x)=
\left\{
\begin{array}{ll}
x,&\mbox{ if   } x<a_\epsilon\\
x+\lambda \frac{(x-a_\epsilon)^4}{(b_\epsilon-a_\epsilon)^4},& \mbox{ if   }x\geq a_\epsilon.
\end{array}
\right.
$$
One can easily see that this function is $C^3$. The perturbation of
the map $f$ will have the form $p_{\epsilon_0,\lambda}\circ f$ for some
sufficiently small $\epsilon_0$ given by the following lemma:

\specialnumber{4.8}\proclaim{Lemma}
  There exist $\lambda_0>0$ and $\epsilon_0$ {\rm (}\/depending on $f${\rm )} such
  that the maps $f$ and $p_{\epsilon_0,\lambda_0}\circ f$ are not
  conjugate and there exists an analytic family of polynomial\/{\rm -}\/like
  maps $F_{\epsilon_0,\lambda}:B_{\epsilon_0,\lambda}\to
  A_{\epsilon_0}$ induced by $p_{\epsilon_0,\lambda}\circ f${\rm ,} where
  $\lambda\in [0,\lambda_0]${\rm ,} $F_{\epsilon_0,0}=F_{\epsilon_0}${\rm ,}
  $B_{\epsilon_0,0}=B_{\epsilon_0}$.
\endproclaim

Before giving a proof of this simple lemma let us notice that though
the map $p_{\epsilon_0}^\lambda\circ f$ is only $C^3$ and not analytic, it
can induce a polynomial-like map because the perturbation is not
analytic just at one point whose forward orbit never comes inside of
$A_{\epsilon_0}$.
\vglue6pt 
\begin{plm}
First of all we can extend the function $p_\epsilon$ to the complex plain by
the following formula:
$$
p_{\epsilon,\lambda}(z)=
\left\{
\begin{array}{ll}
  z,&\mbox{ if   } \Re(z)<a_\epsilon\\
  z+\lambda \frac{(z-a_\epsilon)^4}{(b_\epsilon-a_\epsilon)^4},& 
  \mbox{ if   }\Re(z)\geq a_\epsilon
\end{array}
\right.
$$
This function is discontinuous along the line $\Re(z)=a_\epsilon$.

Fix small $\lambda_0>0$.
Consider a polynomial-like map $F_\epsilon$ and let us see what
happens to it when we perturb the map $f$. 

Due to Theorem~\ref{p_like}, we know that the interval
$(a_\epsilon,b_\epsilon)$ is disjoint from $A_\epsilon$ and that if
$F_\epsilon|_{B_\epsilon^x}=f^n$, then $f^i(x)\notin A_\epsilon^c$ for
$i=1,\dots,n-1$. This implies that if we perturb $f$ by
$p_{\epsilon,\lambda}$, then this will not affect the map $F_\epsilon$
outside of $A_\epsilon\setminus A_\epsilon^c$. Let
$F_{\epsilon,\lambda}|_{B_\epsilon\setminus A_\epsilon^c}=F_\epsilon$.

Again due to Theorem~\ref{p_like} if $x\in B_\epsilon\cap
A_\epsilon^c$, then the size of $f(B_\epsilon^x)$ is very small
compared to $|b_\epsilon-a_\epsilon|$. Hence if $\epsilon$ is small
enough, we have 
$$
f^{-1}\circ {p_{\epsilon,\lambda}}^{-1}\circ
f(B_\epsilon^x) \subset A_\epsilon^c
$$
 for any $x\in B_\epsilon\cap
A_\epsilon^c$, where $0\leq \lambda\leq \lambda_0$.  Let
$$
B_{\epsilon,\lambda}=
\left(B_\epsilon\setminus A_\epsilon^c \right)
\bigcup 
\left(f^{-1}\circ
p_{\epsilon,\lambda}^{-1} \circ f(B_\epsilon\cap A_\epsilon^c)
\right).
$$
 As we have seen, $B_{\epsilon,\lambda}\subset
A_\epsilon$ for $0\leq\lambda\leq\lambda_0$. Finally let
$$
  F_{\epsilon,\lambda}(x)=
  f^{n-1}\circ p_{\epsilon,\lambda} \circ f(x),
$$
where $x\in B_{\epsilon,\lambda}$ and $ n $ is such
that 
$$F_\epsilon|_{B_\epsilon^{(f^{-1}\circ p_{\epsilon,\lambda}\circ
    f(x))}}=f^n.$$
Notice that if $x\notin A_\epsilon^c$, then
$F_{\epsilon,\lambda}(x)=F_\epsilon$.

Decreasing $\epsilon$ if necessary we can get the following:
$f(c)\notin f(B_{\epsilon,\lambda_0})$. Indeed, we know that the ratio
$\frac{|b_\epsilon-f(c)|}{|f(c)-a_\epsilon|}$ can be made arbitrarily
small by decreasing $\epsilon$, so that $p_{\epsilon,\lambda_0}\circ
f(c)\notin f(B_\epsilon)$. Thus $F_\epsilon$ and
$F_{\epsilon,\lambda_0}$ cannot be conjugate.
\end{plm}

Notice that the perturbation $p_{\epsilon_0,\lambda_0}\circ f$ of the
map $f$ is large even in the $C^1$ topology.

The family of polynomial-like maps $F_{\epsilon_0,\lambda}$ is not
trivial: $F_{\epsilon_0,0}$ and $F_{\epsilon_0,\lambda_0}$ are not
conjugate. To this family we can apply the results of
Section~4.2 and conclude that there is
$\lambda_1\in (0, \lambda_0)$ such that the maps $F_{\epsilon_0,0}$ and
$F_{\epsilon_0,\lambda}$ are not conjugate for any
$\lambda\in (0,\lambda_1)$.  Hence, the maps $f$ and $f_\lambda$ are not
conjugate as well, where $f_\lambda=p_{\epsilon_0,\lambda}\circ f$ and
$0<\lambda<\lambda_1$.

We already know that the map $f$ has a $C^3$-neighborhood consisting
of regular maps. Let us denote this neighborhood by $U$.  Taking a
smaller neighborhood if necessary we can assume that $U$ is convex.
Take $\lambda_2<\lambda_1$ so small that the maps $f_\lambda$ belong
to $U$ for $0<\lambda\leq\lambda_2$.  Approximate this map
$f_{\lambda_2}$ by some analytic map $g$ in such a way that the map
$g$ also belongs to $U$ and the maps $g$ and $f$ are not conjugate.
Notice that all the  maps of the  family $f_\lambda$ have a critical point
which does not depend on $\lambda$ and the map $g$ can be chosen in
such a way that the critical points of $f$ and $g$ coincide. Let
$g_\lambda=\lambda g+(1-\lambda)f$, $\lambda\in [0,1]$. Then
$g_\lambda$ is an analytic nontrivial family of analytic unimodal
regular maps with nondegenerate critical point.
Theorem~\ref{real_anal} implies that for small $\lambda$ the maps $f$
and $g_\lambda$ are not conjugate. It is also clear that $f$ and
$g_\lambda$ are close in the $C^\omega(\Delta)$ topology for small
$\lambda$.  \hfill\qed

 \section{Appendix}
\label{sec:appendix}

5.1. {\it Quasiconformal homeomorphisms}.
In this section we will give a short overview of definitions and
results connected with quasiconformal maps. For the details the
reader can consult books \cite{Ahl}, \cite{LV}.

There are many different, equivalent definitions of the quasiconformal
(q.c.) homeomorphism. We will use the following:
\numbereddemo{Definition}
  Let $U\subseteq \bar\Cp$ be a domain in the complex plane. The map
  $h:U\to h(U)$ is called a {\it quasiconformal} homeomorphism if
  \begin{itemize}
  \item $h$ is an orientation preserving homeomorphism between the
    domains $U$ and $h(U)$;
  \item The real part $\Re(h)$ and the imaginary part $\Im(h)$ of $h$
    are absolutely continuous on almost all verticals and almost all
    horizontals in the sense of Lebesgue;
  \item There exists a constant $k<1$ such that for 
    $$
    \mu_h(z)=\frac{d_{\bar z}f(z)}{d_zf(z)}
    $$
    one has 
    $$
    |\mu_h(z)|<k
    $$
    for almost all $z\in U$ where $d_{\bar z}h=\frac{dh}{d\bar z}$ and
    $d_zh=\frac{dh}{dz}$.
  \end{itemize}
\enddemo

The function $\mu_h$ is called the {\it Beltrami coefficient} of a
q.c.~homeomorphism~$h$.

To the Beltrami coefficient $\mu$ one can associate a field of
infinitesimal ellipses. The eccentricities of these ellipses are given
by $\frac{1+|\mu(z)|}{1-|\mu(z)|}$ and the  directions of the major
axes are given by $\sqrt{\mu(z)}$.

If $f$ is a holomorphic map, we can pull back this field of ellipses
even if $f$ is not injective. This pullback we will denote as
$f^*\mu$ which is equal to
$$
(f^*\mu)(z)=\mu(f(z)) \frac{\overline{d_zf(z)}}{d_zf(z)}.
$$

Here is a list of theorems to be used later on.

\advance\theoremcount by -1

\proclaimtitle{measurable Riemann mapping theorem} \label{mRm}
\proclaim{Theorem}
  Let $\mu:\Cp \to \Cp$ be a measurable function such that $|\mu|<k<1$
  almost everywhere. Then there exists a unique {\rm q.c.}~homeomorphism
  $h:\bar \Cp\to \bar\Cp$ whose Beltrami coefficient is $\mu$ and which is
  normalized such that $h(0)=0$, $h(1)=1$ and $h(\infty)=\infty$. 
\endproclaim

\proclaimtitle{Ahlfors-Bers theorem} 
\proclaim{Theorem} \label{ahl_be}
  Let $\Lambda\subset \Cp^n$ be an open set and
  $\mu:\Cp\times \Lambda\to \Cp$ be a measurable function satisfying\/{\rm :}
  \begin{itemize}
  \item $|\mu(z,\lambda)|<k<1$ for all $\lambda\in \Lambda$ and for
    almost all $z\in \Cp$;
  \item The map $\lambda\mapsto \mu(z,\lambda)$ is holomorphic in
    $\lambda$ for almost all $z\in \Cp$.
  \end{itemize}
  Then there exists a unique function $H:\Cp\times \Lambda\to \Cp$ such
  that 
  \begin{itemize}
  \item $H(0,\lambda)=0${\rm ,} $H(1,\lambda)=1${\rm ,}
    $H(\infty,\lambda)=\infty${\rm ;}
  \item For fixed $\lambda\in \Lambda$ the map $z\mapsto F(z,\lambda)$
    is a {\rm q.c}.~homeomorphism whose Beltrami coefficient is
    $\mu(\cdot,\lambda)${\rm ;}
  \item The map $\lambda\mapsto F(z,\lambda)$ is holomorphic for
    almost every $z$.
  \end{itemize}
\endproclaim 

\vglue-12pt
The first version of the next theorem appeared in \cite{MSS} and
after it was generalized several times: \cite{BR}, \cite{Slo}.

\proclaimtitle{$\lambda$-lemma}
\proclaim{Theorem} \label{la-lemma}
  Let $Z\subset \bar\Cp$ be a set{\rm ,} $D$ be an open unit disk in the
  complex plane and let $h:Z \times D\to \bar\Cp$ satisfy the following
  conditions\/{\rm :}
  \begin{itemize}
  \item $h(z,0)=z$ for any $z\in Z${\rm ;}
  \item For fixed $z\in Z$ the function $\lambda\mapsto h(z,\lambda)$
    is holomorphic for $\lambda\in D${\rm ;}
  \item For fixed $\lambda\in D$ the map $z\mapsto h(z,\lambda)$ is
    injective for all $z\in Z$.
  \end{itemize}
  Then there exists $H:\bar\Cp \times D \to \bar\Cp$ such that
  \begin{itemize}
  \item $H(z,\lambda)=h(z,\lambda)$ for $\lambda\in D$ and $z\in Z${\rm ;}
  \item $H(z,0)=z$ for $z\in \bar \Cp${\rm ;}
  \item For fixed $z\in \bar\Cp$ the function $\lambda\mapsto
    H(z,\lambda)$ is holomorphic for $\lambda\in D${\rm ;}
  \item For fixed $\lambda\in D$ the map $z\mapsto H(z,\lambda)$ is a
    {\rm q.c.} homeomorphism\/{\rm ;}
  \item For almost every $z\in \bar\Cp$ the Beltrami coefficient of
    $H$ depends holomorphically on $\lambda$.
  \end{itemize}

\endproclaim 
\vglue-12pt
Since the Beltrami coefficient of a q.c.~homeomorphism is not defined
everywhere we have to clarify the last item in the previous theorem.
We say that the Beltrami coefficient depends holomorphically on
$\lambda$ for almost every $z$ if there is a function $\mu(z,\lambda)$
such that for almost every $z$ the function $\lambda\mapsto
\mu(z,\lambda)$ is holomorphic and for fixed $\lambda$ the equality
$\mu(z,\lambda)=\mu_{H(\lambda,\cdot)}(z)$ holds almost everywhere.

\proclaimtitle{Compactness of the set of q.c.~homeomorphisms}
 \proclaim{Theorem}  If $H$ is a family of {\rm q.c}.~homeomorphisms of $\bar\Cp$ whose
  Beltrami coefficients are uniformly bounded by a constant $k<1${\rm ,} then
  any sequence in $H$ has a subsequence which converges uniformly and
  the limit either a constant or a {\rm q.c.}~homeomorphism whose Beltrami
  coefficient is bounded by $k$.  
\endproclaim 

\proclaim{Theorem}
  If $f$ is holomorphic{\rm ,} then  $\mu_{f\circ h}=\mu_h$ and $\mu_{h\circ
  f}(z)=\mu_h(f(z))\frac{\overline{d_zf(z)}}{d_zf(z)}$.
\endproclaim 

The real counterpart of q.c.~homeomorphisms are quasisymmetric
homeomorphisms of the real line. 

\demo{Definition {\rm 5.2}}
  The homeomorphism $h:\Rp\to \Rp$ is called {\it quasisymmetric} if
  there is a constant $C>0$ such that for any three points
  $x_{-1}<x_0<x_1$ such that $x_0-x_{-1}=x_1-x_0$ the following
  inequality holds:
  $$
  C^{-1}<\frac{|h(x_1)-h(x_0)|}{|h(x_0)-h(x_{-1})|}<C.
  $$
\enddemo

The following theorem describes relations between quasiconformal and
quasisymmetric homeomorphisms:

\proclaim{Theorem}
  Let $h^c$ be a quasiconformal homeomorphism of the complex plane such
  that its restriction $h^r$ to the real line is a real function. Then
  this restriction is a quasisymmetric homeomorphism.

  If $h^r$ is a quasisymmetric homeomorphism of the real line{\rm ,} then
  there is a quasiconformal homeomorphism $h^c:\Cp\to \Cp$ such that
  the restriction of $h^c$ to the real line is $h^r$.
\endproclaim

5.2. {\it The straightening theorem and geodesic neighborhoods}.
One of the important applications of the measurable Riemann mapping
theorem to holomorphic dynamical systems is the straightening theorem.
Let $f:B\to A$ be a holomorphic proper 2-to-1 map where $B$ and $A$
are simply connected domains and $A$ contains the closure of $B$. Such
a map is called {\it quadratic-like}. Let $J(f)=\{z\in \Cp: f^i(z)\in
U \mbox{ for all } i\leq 0\}$.  This set is called the {\it filled
  Julia set } of the quadratic-like map $f$. Douady and Hubbard proved the following result:

\proclaimtitle{The straightening theorem \cite{DH}} 
\proclaim{Theorem}\label{straigh} 
  Let $f:B\to A$ be a quadratic\/{\rm -}\/like map and $d$ be the degree of $f$.
  Then there exists a quadratic map $p${\rm ,} a neighborhood $U$ of $J(f)$
  such that $f:U\to f(U)$ is a quadratic\/{\rm -}\/like map and there is a
  {\rm q.c.}~homeomorphism $h:f(U)\to p(h(U))$ which conjugates $f|_U$ and
  $p|_{h(U)}$.
\endproclaim 



Let $I$ be some interval on the real line. $\Cp_I$ will denote the
domain $\Cp\setminus(\Rp\setminus I)$. Consider the Poincar\'e metric
on the domain $\Cp_I$. It is clear that $I$ is a geodesic in this
metric. Denote the set of points whose distance in this metric to the
interval $I$ is less than $l$ as $\tilde D_l(I)$. 

Consider two circles $S^-$ and $S^+$ centered at the points $a^-$ and
$a^+$ such that these points are symmetric with respect to the real
line, and let these circles pass through the boundary points of the
interval \pagebreak $I$ and intersect the real line at the angle
$\phi<\frac{\pi}{2}$. Denote the intersection of the disks delimited
by these circles as $D_\phi(I)$ and the union of these disks as
$D_{\pi-\phi}(I)$. So, $D_\phi(I)$ is a lens as shown in
Figure~1.

One can check that the domain $\tilde D_l(I)$ coincides with $D_\phi(I)$ for
$l=\ln\tan(\frac{ \pi}{4} +\frac{\phi}{4})$. (See \cite{dMvS}.) 

If $g$ is a univalent map of the domain $\Cp_I$, then it contracts the
Poincar\'e metric. So we have the following lemma:

\specialnumber{5.1}\proclaim{Lemma}
  Let $g:\Cp_I\to \Cp_{g(I)}$ be a univalent map and let $g(I) \subset
  \Rp$. Then for any interval $J\subseteq I$ and any $\phi${\rm ,}
  $$
  g(D_\phi(J))\subseteq D_\phi(g(J)).
  $$
\endproclaim

Obviously, if the interval $I$ consists of positive real numbers,
then the square root map is univalent on $\Cp_I$ and we can apply the
previous lemma. Another case when we can use it, is a case of the
Epstein class.

\demo{Definition {\rm 5.3}} 
  A map $f$ belongs to the Epstein class  if it is real
  analytic and any inverse branch $f^{-1}:I\to \Rp$ can be univalently
  extended to the domain $\Cp_I$; i.e.,  if $J$ is an interval
  of the monotonicity of $f$ and $I=f(J)$, then the map $f^{-1}|_I$
  can be holomorphically extended and the extended map
  $f^{-1}:\Cp_I\to\Cp_J$ is univalent.
\enddemo

If an analytic map does not belong to the Epstein class, whenever the
size of $D_\phi(I)$ is small compared to the size of the  extension of
$f^{-1}$ to the complex plain, one can give an estimate of the shape
of the pullback of $D_\phi(I)$. More precisely, the following lemma
holds:

\proclaimtitle{[dFdM,  Lemma~2.4]}
\specialnumber{5.2 }\proclaim{Lemma}
  \label{lm:fm}
  There exists a universal constant $\tau_3>1$ such that for any
  small $a>0$ there exists $\theta(a)\in (0,\pi)$ satisfying
  $\theta(a)\to \pi$ and $a/(\pi-\theta(a))\to 0$ as $a\to 0$ such
  that the following holds. Let $F:D\to \Cp${\rm ,} where $D$ is a unit
  disk{\rm ,} be univalent and symmetric with respect to the real line{\rm ,} and
  assume that $F(0)=0${\rm ,} $F(a)=a$. Then for all $\phi\in (0,\theta(a))${\rm ,}
  $$
  F(D_\phi([0,a]))\subset D_{(1+a^{\tau_3})\phi}([0,a]).
  $$
\endproclaim

 \newcommand{\ef}{{\cal E}}
\newcommand{\rf}{{\cal R}}
\newcommand{\cf}{{\cal C}}
\newcommand{\wf}{{\cal W}}

5.3. {\it Construction of the holomorphic box mapping}.
Following a suggestion of the referee  we include an outline of the proof of
Theorem~\ref{thr:lvs} here. This theorem was proved in
\cite[Th.~C]{LvS} in the case of maps of the form $x\mapsto x^l+c$
where $l$ is even and $c$ is real. To generalize the result of
\cite{LvS} we will follow the proof given in Section~14 of
\cite{LvS}. We will also use the notation of that paper (though the
author of the present paper thinks that it is slightly illogical) even
if it is different from what we have used above.  Though \pagebreak we will not
give proofs of lemmas if they are identical to \cite{LvS} we will
try to keep the exposition self-contained. In what follows we will
assume that $f$ is nonrenormalizable since the renormalizable case appears in \cite[Th.~11.1]{LvS}.

Given a unimodal map $f$ we say that  $g\in \ef(T^0)$  if $T^0$ is a
nice symmetrical interval around the critical point and $g:\cup_i
T^1_i\to T^0$ where $\cup_i T^1_i$ is a collection of disjoint
subintervals of $T_0$. Moreover, the following properties are
satisfied:
\begin{itemize}
\item If $i\neq 0$, the map $g:T^1_i\to T^0$ is a diffeomorphism onto
  $T^0$ of the form $f^{j(i)}$;
\item Denoting $T^1_0$ by $T^1$ we have that $g:T^1\to T^0$ is a
  unimodal map of the form $f^j$, $g(\partial T^1)\in T^0$ and the range of
  the map $f^{j-1}:f(T^1)\to T^0$ can be extended to $T^0$;
\item All iterates of the critical point under $g$ are in $\cup_i T^1_i$.
\end{itemize}

Next we say that $g\in \ef(T^0,T^{-1})$ if $T^{-1}$ is a nice
symmetrical interval containing $T^0$, $g\in \ef(T^0)$ and the range
of the map $f^{j(i)-1}:f(T^1_i)\to T^0$ can be extended to $T^{-1}$
for all $i$.

We can define low, high and center returns for maps in $\ef(T^0)$ in
the same way we did it for first entry maps in Section~1.6.

Now we introduce a renormalization operator $\rf$ for maps in
$\ef(T^0)$. Notice that this operator is different from the one
used above.

First, we define $\rf g$ in the case when $g$ is a low return. In this
case $\rf g$ will be in $\ef(T^0)$.  

Let $g$ be a low return. For any point $x \in \cup T^1_i$ we define
$s(x)$ as a minimal nonnegative integer such that $g^{s(x)}(x) \notin T^1$
(thus for $x \notin T^1$ we have $s(x)=0$). Then we define the
intermediate renormalizations $\hat\cf g$ by $\hat\cf
g(x)=g^{s(x)+1}(x)$ and $\cf g$ by $\cf g(x)=g(x)$ if $x\notin T^1$
and $\cf g(x)=g^{s(x)}(x)$ if $x\in T^1$ (the definition of $\cf g$ is given just to keep the same notation as
in
\cite{LvS}; we will not use it).

\specialnumber{5.3}\proclaim{Lemma}
  If $g\in \ef(T^0)${\rm ,} then the map $\hat\cf g$ is in $\ef(T^0)$ as
  well.
\endproclaim

\begin{plm}
It is easy to see that any noncentral branch of $\hat\cf g$ is a
diffeomorphism onto $T^0$. Let $\hat T^1\subset T^1$ be a central domain of
$\hat \cf g$ and $V$ be a domain of $\hat\cf g$ such that $g(c)\in
V$. Then $\hat \cf g|_{\hat T^1}=\hat\cf g|_V \circ g|_{\hat
  T^1}$. However, since the range of the map $f^{j-1}:\hat T^1\to
T^0$ can be extended to $T^0$, where $g|_{T^1}=f^j$,  the
interval $V$ is contained in this range. Now using the fact that
$\hat \cf g|_V$ is a diffeomorphism onto $T^0$ we obtain that the
range of the map $\hat \cf g|_{\hat T^1}$ can be extended to
$T^0$.
\end{plm}
\vglue-9pt
If $\hat\cf g $ is a low return again, we can define $\hat\cf^2 g$
(for the second intermediate renormalization the function $s$ has to
be defined with respect to $\hat T^1$) and so on. Let $\hat T^i$  \pagebreak be a
sequence of central domains of $\hat\cf^i g$ and let $\tilde T^1$ be
the central domain of $\cf g$. Let $\tilde s$ be
minimal nonnegative number such that $\hat\cf^{\tilde s} g(\hat
T^{\tilde s}) \cap \tilde T^1 \neq\emptyset$. Then the renormalization of
$g$ is  $\rf g=\hat\cf^{\tilde s} g$. As a consequence of the
previous lemma we obtain that $\rf g\in\ef(T^0)$.

Now let $g$ be a high return. Let $x$ be an orientation preserving
fixed point of $g|_{T^1}$ and $z_1$ be a boundary point of $T^1$ such
that $x$ is between $c$ and $z_1$. Take preimages $z_2,z_3,\ldots$ of $z_1$
along the branch $g|_{[z_1,c]}$. Let $U_k$ be an interval with
boundary points $z_k$ and the point symmetrical to $z_k$ and choose
$k\geq 0$ minimal such that $g(U_k)\supset U_k$. The map $f$ is not
renormalizable, hence $k$ exists. Denote $U_k$ by $V^1$. For $x\notin
V^1$ define $\tilde \wf g(x)$ by $\tilde \wf g(x)=g^j(x)$ where $j$ is
minimal such that $x\notin U_j$. For $x\in V^1$ let $\tilde \wf g$ be
the first return map of $g$ to $V^1$. Finally, let $\wf g=\tilde \wf
g|_{V^1}$. 

In \cite[Lemma~14.1]{LvS} it is proved that if $g\in \ef(T^0)$ and $g$
is a low return, then $\wf g\in \ef(V^1,T^0)$.

\specialnumber{5.4}\proclaim{Lemma}
  Suppose that $g\in\ef(T^0)$ is a first return map of $f$ to $T^0$ and
  that $\hat\cf^i g$ is a low return for $i=0,\ldots,m-1$. Let $U$ be a
  domain of $\hat\cf^m g$ and let $\hat\cf^m g|_U=f^n$. Then the
  orbit $f(U),f^2(U),\ldots,f^m(U)$ has intersection multiplicity at most
  $m+1$.
\endproclaim

Here we say that a collection of intervals has {\it intersection
  multiplicity} $k$ if any point is covered by not more than $k$
intervals from this collection.

This lemma can be proved  easily  by induction.

\specialnumber{5.5}\proclaim{Lemma}\label{lm:fret}
  Suppose that $g\in\ef(T^0)$ is a first return map of $f$ to $T^0${\rm ,}
  that $\hat\cf^i g$ is a low return for $i=0,\ldots,m-1$ and let $\hat
  T^i$ be a central domain of $\hat\cf^i g$. Then the first return
  map of $\hat\cf^m g$ to $T^m$ coincides with the first return map
  of $f$ to $T^m$.
  
  Moreover{\rm ,} if $\hat\cf^m g$ is a high return{\rm ,} then $\wf\circ
  \hat\cf^m g \in\ef(V^1,T^0)$ is a first return map of $f$ to the
  interval $V^1$.
\endproclaim

\begin{plm} 
  We will prove by induction with respect to $m$ that the first return
  map  of $\hat\cf^m g$ to any nice interval $U$ contained in
  $\hat T^m$ is the first return map of $f$ to~$U$. 
  
  Let $x\in U$.  Let $R$ be the first return map of $\hat\cf^{m-1}
  g$ to $U$. By the induction assumption $R$ coincides with the first
  return map of $f$ to $U$.  Let $R(x)=(\hat\cf^{m-1} g)^n(x)$. Then
  $(\hat\cf^{m-1} g)^{n-1}(x)\notin\hat T^{m-1}$ because by the
  construction of $\hat \cf^{m-1} g$ we have $\hat\cf^{m-1} g(\hat T^{m-1})\cap
  U=\emptyset$. Thus $R(x)$ can be written as
  \begin{eqnarray*}
    R(x)&=&(\hat\cf^{m-1} g)^n(x) \\ 
    &=&
    (\hat\cf^{m-1} g|_{T^0\setminus\hat T^{m-1}}\circ
    (\hat\cf^{m-1} g|_{\hat T^{m-1}})^{s_l})\\
&&\circ\cdots\circ(\hat\cf^{m-1}
    g|_{T^0\setminus\hat T^{m-1}} \circ (\hat\cf^{m-1} g|_{\hat
      T^{m-1}})^{s_1})(x)\\
    &=&(\hat\cf^m g)^l(x),
  \end{eqnarray*}
  where $s_i\geq 0$, $i=1,\ldots,l$.
  Therefore, $R$ is the first return map of $\hat\cf^n g$ to $U$.

  The case of the high return can be treated in the same way.
\end{plm}

\specialnumber{5.6}\proclaim{Lemma}\label{lm:real-b2}
  Let $f$ be a $C^3$ unimodal map with a nondegenerate recurrent
  critical point. For any $\epsilon>0$ there exists $\tau<1$ such that if
  $T^{-1}$ is a sufficiently small nice interval{\rm ,} $g\in\ef(T^0,T^{-1})${\rm ,}
  $T^{-1}$ is an $\epsilon$\/{\rm -}\/scaled neighborhood of $T^0${\rm ,} $\hat\cf^i g$ is a low
  return for $i=0,\ldots,m-1${\rm ,} $T^i$ is a central domain of $\rf^i g${\rm ,} then
  $$
  \frac{|T^i|}{|T^0|}<\tau^i
  $$
  where $i=1,\ldots,m$.
\endproclaim

\begin{plm}
  The standard cross-ratio estimate yields the fact that for any $\epsilon>0$ there
  is $K<1$ such that if $T^1_j$ is a domain of $g$, then
  $\frac{|T^1_j|}{|T^0|}<K$. Applying the standard cross-ratio
  estimate once again we obtain the required inequality. \phantom{whatever}\hfill
\end{plm}

The next three lemmas are a version of Lemma~14.4 of \cite{LvS}
broken into three parts and adapted to our case.

\specialnumber{5.7}\proclaim{Lemma}\label{lm:14.4-1}
  Let $f$ be a $C^3$ unimodal map with a nondegenerate recurrent
  critical point. For any $\epsilon>0${\rm ,} $\tau>0$ there exists $N$ such that if
  $T^{-1}$ is a sufficiently small nice interval which is an $\epsilon$\/{\rm -}\/scaled
  neighborhood of $T^0${\rm ,} $g\in\ef(T^0,T^{-1})$ is a first return of $f$
  to $T^0${\rm ,} $\hat \cf^i g$ is a low return for $i=0,\ldots,N${\rm ,} $\hat T^N$ is a
  central domain of $\hat\cf^N g${\rm ,} $R$ is the first return map of $f$ to
  $\hat T^N${\rm ,} $U$ is a central domain of $R$, $R|_U=f^j${\rm ,} then the range of
  the map $f^{j-1}:f(U)\to T^N$ can be extended to $T^0$. Moreover{\rm ,} if
  $W$ is a connected component of the preimage $f^{-j+1}(T^N)$
  containing $f(U)${\rm ,} then
  $$
  \frac{|W\setminus f(U)|}{|f(T^N)|}<\frac 12.
  $$
\endproclaim

\begin{plm}
  First,  we notice that $f^j(\partial T^N)$ is not in the interior of
  $T^0$. Indeed, due to Lemma~\ref{lm:fret},   $R$ is a first return map
  of $\hat\cf^N g$ to $\hat T^N$, so that  $f^j=f^{j_1}\circ\hat\cf^N g$ for
  some $j_1\geq 0$.  Now, $\hat\cf^N(\partial T^N)\in \partial T^0$ and $T^0$ is a nice
  interval; hence $f^j(\partial T^N)\notin\mathop{\mbox{int}}T^0$.

  Next, by standard arguments (e.g. see Lemma~\ref{exten_centr}) we
  obtain the extension of the range of $f^{j-1}:f(U)\to T^N$. 
  
  The required inequality can be obtained by the same estimate as in\break
  Lemma~\ref{real_b}. Note that we can use this estimate because the
  preimage of one of the boundary points of $T^0$ by $f^{-j+1}$ is in
  the closure of $f(T^N)$.
\end{plm}

\specialnumber{5.8}\proclaim{Lemma}\label{lm:14.4-2}
  Let $f$ be a $C^3$ unimodal map with a nondegenerate recurrent
  critical point. For any $N$ and $\epsilon>0$ there is $\sigma\in(0,1)$ such that
  if $T^{-1}$ is a sufficiently small nice interval{\rm ,} $g\in
  \ef(T^0,T^{-1})$ where $T^{-1}$ is an $\epsilon$\/{\rm -}\/scaled neighborhood of
  $T^0${\rm ,} $\rf^i g$ is a low return for $i=0,\ldots,k-1$ where $k<N${\rm ,}
  $\rf^k$ is a high return{\rm ,} $T^{k+1}$ is a central domain of
  $\rf^k g${\rm ,} $\rf^k g|_{T^{k+1}}=f^j${\rm ,} $A$ is a connected
  component of the preimage $f^{-j+1}(T^0)$ containing $f(T^{k+1})${\rm ,}
  and
  $$
  \frac{|T^0|}{|V^1|}<(1-\sigma)^{-1} \frac{|T^{-1}|}{|T^{0}|},
  $$
  then
  $$
  \frac{|A\setminus f(T^{k+1})|}{|f(T^0)|}<1-\sigma.
  $$
 \endproclaim

\begin{plm}
  The proof of this lemma is nearly identical to the proof of the second
  assertion of Lemma~14.4 of \cite{LvS} with slight modifications. 

  All combinatorial properties of unimodal maps used in that proof
  obviously hold in our case. In particular, Lemma~14.2 holds. The
  estimates of that proof also hold with some spoiling factors
  close to one. Namely, the  first first spoiling factor appears in
  inequality~14.2, which in our case would look like:
  $$
  \frac{|R|}{|I|}> C \mu_{i+1} \frac{|f(R)|}{|f(I)|}
  $$
  where $C$ is a constant close to $1$ if $T^0$ is small. In the
  same fashion such spoiling factors appear in other inequalities
  there. They will start accumulating as we have more and more low
  returns. Thus, inequality~14.7 will have the form
  $$
  \frac{|R'|}{|I'|}>
  C^r (1-\epsilon(\sigma)) \left(\frac{|T^{k+1}|}{|T^0|}\right)^2
  \left(\frac{|T^0|}{|T^1|}\right)^\tau.
  $$
  Here in our notation $R'=f(T^{k+1})$ and $I'=A\setminus f(T^{k+1})$.

  The number of low returns is bounded by $N$, hence $r\leq N$ is also
  bounded. If the interval $T^0$ is small enough, the constant $C$ can
  be made as close to $1$ as we want. Therefore, we get the same
  estimate~14.8
  $$
  \frac{|R'|}{|I'|}>
  \kappa \left(\frac{|T^{k+1}|}{|T^0|}\right)^2
  $$
  where $\kappa>1$ is some constant. The rest is the same as in
  \cite{LvS}.
\end{plm}

\specialnumber{5.9}\proclaim{Lemma}\label{lm:14.4-3}
  Let $f$ be a $C^3$ unimodal map with a nondegenerate recurrent
  critical point. There is $\epsilon>0$ such that if $g\in \ef(T^0,T^{-1})$
  is a first return map of $f$ to a sufficiently small interval $T^0${\rm ,}
  $T^{-1}$ is an $\epsilon$\/{\rm -}\/scaled neighborhood of $T^0${\rm ,} $g|_{T^1}=f^j${\rm ,} $W$
  is a connected component of the preimage $f^{-j+1}(T^0)$ containing
  $f(T^1)${\rm ,} then
  $$
  \frac{|W\setminus f(T^1)|}{|f(T^0)|}<\frac 12.
  $$
\endproclaim

This proof is identical to the proof of Lemma~\ref{lm:14.4-1}.

\demo{Proof of Theorem {\rm 2.3}}  Let $f$ be a real-analytic nonrenormalizable
unimodal map with nondegenerate recurrent critical point. After some
analytic change of the coordinate we can assume that $f=\hat f(x^2)$
where $\hat f$ is a real-analytic diffeomorphism. Let $\Omega$ be a complex
neighborhood of the image of $f$ such that $\hat f^{-1}$ is univalent
on $\Omega$.

We know that there is a sequence of pairs of nice intervals
$\{T^0_i,T^{-1}_i\}$ whose lengths tend to zero and such that the first
return map of $f$ to $T^0_i$ is in $\ef(T^0_i,T^{-1}_i)$. Moreover,
there exists a constant $\epsilon>0$ such that for all $i$ the interval
$T^{-1}_i$ is an $\epsilon$-scaled neighborhood of $T^0_i$ (see
\cite{Mar} or \cite{Koz}).

For this $\epsilon$, Lemma~\ref{lm:14.4-1} gives $N$. There is also $\delta>0$ such
that if $U$ is a domain of the first return map to $T^0_i$, then
$T^0_i$ is a $\delta$-scaled neighborhood of $U$. Fix some angle $\phi_0$
slightly less than $\frac \pi2$. Then there is $\phi_1\in(\phi_0,\frac \pi2)$
such that $D_{\phi_1}(U)\subset D_{\phi_0}(T)$ if $T$ is a $\delta$-scaled neighborhood
of $U$. Moreover, the modulus of $D_{\phi_0}(T)\setminus D_{\phi_1}(T)$ is bounded
away from zero by some constant which depends only on $\delta$. 

Take a pair $\{T^0,T^{-1}\}$ from the sequence with such small intervals
that Lemmas~\ref{lm:real-b2}, \ref{lm:14.4-1}, \ref{lm:14.4-2} and
\ref{lm:14.4-3} start to work. Moreover, let $T^0$ be so small that if
$f^n(U)\subset T^0$, $f^n|_U$ is a diffeomorphism, the intersection
multiplicity of the orbit $f(U),\ldots,f^n(U)$ is at most $N$, then
$\sum_{i=1}^n |f^i(U)|^{\tau_3}<\log(\phi_1/\phi_0)$ and intervals $f^i(U)$
are small compared with the distance to $\partial \Omega$ as Lemma~\ref{lm:fm}
requires. Here the constant $\tau_3>1$ is given by
Lemma~\ref{lm:fm}. Such a  $T^0$ exists because of the absence of
wandering domains; see also Lemma~5.2 in \cite{Koz}. The last
inequality implies $\prod_{i=1}^n(1+|f^i(U)|^{\tau_3})<\phi_1/\phi_0$.

Let $g$ be the first return map to $T^0$. If $g$ is a low return as
well as $\hat\cf^i g$ for $i=1,\ldots,N$, then let $R$ be the first return
to $\hat T^N$ where $\hat T^N$ is a central domain of $\hat\cf^N g$.
Let $U$ be any noncentral domain of $R$ and $R|_U=f^j$. Then the
orbit of $U$ is disjoint and $f^{-j}(D_{\phi_0}(T^N)\subset
D_{\prod_{i=1}^j(1+|f^i(U)|^{\tau_3})\phi_0}(U)\subset D_{\phi_1}(U)\subset D_{\phi_0}(T^N)$.
For the central domain $U$ we have $f^{-j+1}(D_{\phi_0}(T^N))\subset
D_{\phi_1}(W)$ where $W$ is as in Lemma~\ref{lm:14.4-1}. Pulling back
$D_{\phi_1}(W)$ by $f$ and using Lemma~\ref{lm:14.4-1} we obtain
$f^{-j}(D_{\phi_0}(T^N))\subset D_{\phi_0}(T^N)$. Notice that all these pullbacks
do not intersect because $\phi_1<\frac \pi2$.

Now let $\rf^i g$ be a low return for $i=0,\ldots,k-1$ and $\rf^k g$ be a high return where $k<N$. By the
construction of $\rf^i g$ we know that
$\rf^k g=\hat\cf^m g$ for some $m$. We can assume that $m<N$;
otherwise we are in the previous case. Suppose we are in the setting
of Lemma~\ref{lm:14.4-2}. Arguing as above we can construct a
holomorphic box mapping whose real trace is $\rf^k g$. Notice that we
can use Lemma~\ref{lm:fm} because orbits of domains of $\rf^k g$ have
intersection multiplicity at most $m+1\leq N$.

If Lemma~\ref{lm:14.4-2} does not apply; i.e., the inequality $
\frac{|T^0|}{|V^1|}\geq(1-\sigma)^{-1} \frac{|T^{-1}|}{|T^{0}|} $ is
satisfied, we can consider the map $g_1=\wf\circ\rf^k g\in\ef(V^1,T^0)$. If
either Lemma~\ref{lm:14.4-1} or Lemma~\ref{lm:14.4-2} applies to
$g_1$, we are done, otherwise we obtain a map $g_2\in\ef(V^2,V^1)$ and
so on. The ratio
$\frac{|V^{i-1}|}{|V^{i}|}>(1-\sigma)^i\frac{|T^{-1}|}{|T^0|}$ tends to
infinity. Thus, sooner or later we will have that the interval
$V^{i-1}$ is an $\epsilon$-scaled neighborhood of $V^i$ where $\epsilon$ is given by
Lemma~\ref{lm:14.4-3}. In this case we can proceed exactly in the same
way as before.  \enddemo

{\it Remark}. We have not used the fact that the critical point is
quadratic. So, one can remove the condition on the nondegeneracy of
the critical point in Theorem~\ref{thr:lvs}. (Note that if $f$ is
real-analytic, the critical point is always nonflat.)

\demo{{\rm 5.4.} Lebesgue measure of the Julia set}
 The following result is proven in \cite[Cor.~2]{Lyu1}:
\enddemo

\proclaimtitle{\cite{Lyu1}}
\proclaim{Theorem}
  \label{thr:lm0}
  Let $F:B\to A$ be a polynomial\/{\rm -}\/like map{\rm ,} where $A$ consists only
  from one connected domain and $\partial B\cap \partial A=\emptyset${\rm ,} and let $F$ be a
  nonremormalizable map. Moreover{\rm ,} suppose that the critical point of
  $F$ is recurrent. Then the Julia set $J=\{x\in B:\, F^n(x)\in B\,
  \forall n\geq 0$ of $F$ has zero Lebesgue measure.
\endproclaim 

Here, $F$   nonrenormalizable means that $F$ does not induce a
quadratic-like map.


This theorem can be easily generalized to arbitrary polynomial
nonrenormalizable maps:

\proclaim{Theorem}
  \label{thr:lm0-general}
  Let $F:B\to A$ be a polynomial\/{\rm -}\/like map{\rm ,} the critical point of $F$ be
  recurrent and let $F$ be a nonremormalizable map. Then the Julia
  set $J$ of $F$ has zero Lebesgue measure.
\endproclaim 

The proof of this statement is very similar to the proofs of similar
statements in \cite{LvS} and \cite{Lyu1}. Following a suggestion of
the referee it is included here.

\begin{plm}
  Consider two cases. First assume that the $\omega$-limit set of $c$ is
  minimal. Let $R:\hat B\to A^c$ be the  first return map to
  $A^c$. The domain $\hat B$ can contain infinitely many connected
  components. However, there are only finitely many connected
  components which contain points of the orbit of the critical
  point. Indeed, since $c$ is recurrent the $\omega$-limit set contains
  the orbit of $c$; since $\omega(c)$ is minimal and does not contain
  points of the boundary of $\hat B$, thus $\omega(c)\subset\hat B$; the set
  $\hat B$ is open and $\omega(c)$ is compact; hence there are finitely
  many components  of $\hat B$ which cover $\omega(c)$.
  
  Let $\tilde B$ be a union of these finitely many components. Then we
  can apply Theorem~\ref{thr:lm0} to the map $R|_{\tilde B}$. So the
  Julia set of $R|_{\tilde B}$ has zero Lebesgue measure, hence the Lebesgue
  measure of $J$ is zero as well.
  
  Now suppose that the $\omega$-limit set of $c$ is nonminimal. Then
  there is a point $a\in \omega(c)$ such that $c\notin \omega(a)$.
  
  Denote $A_0=A$, $A_1=B$ and $A_k=F^{-k}(A)$. The map $F$ is
  nonrenormalizable; hence sizes of domains in $A_k$ shrink to zero.
  Let ${k_0}$ be such that $A_{k_0}^c$ does not contain points from
  the orbit of $a$ and let $k_1$ be such that $A_{k_1}^c$ is compactly
  contained in $A_{k_0}^c$. Consider the first entry map $R:\hat B\to
  \hat A$, where $\hat A=A_{k_1}^c$. It is easy to see that if $x\in
  (\hat B\setminus\hat A)$ then the range of $R$ can be univalently extended
  to $A_{k_0}^c$ (compare Lemma~\ref{exten_centr}). On the other hand,
  we also have the following property: let the map $R:\hat B^x\to \hat
  A$ have a univalent extension $\tilde R:\tilde B\to A_{k_0}^c$; then
  for any $k$ either $A_k^x$ contains $\tilde B$ or $\tilde B$
  contains $A_k^x$.

  Suppose that the Julia set $J$ has positive Lebesgue measure. Let $b$ be
  a density point of $J$ such that $c\in \omega(b)$. Such a point always
  exists because the set of points whose $\omega$-limit sets do not
  contain the critical point has zero Lebesgue measure. Given $k$ let
  $n_k$ be a minimal integer such that $F^{n_k}(b)\in A_k^a$. Such an
  integer $n_k$ exists because $a\in \omega(c)\subset\omega(b)$. Since
  $n_k$ is minimal there is a domain $U_k\ni b$ such that $F^{n_k}$ maps
  $U_k$ onto $A_k^a$ either univalently or as\break 2-to-1. Then there is
  domain $V_k\ni b$ such that $R\circ F^{n_k}$ maps $V_k$ onto $\hat
  A$ univalently or as 2-to-1. The range of this map can be extended
  to $A_{k_0}^c$ and again the extension map is at most 2-to-1. As
  $k\to \infty$ the size of $V_k$ tends to zero. Since $b$ is a density
  point of $J$ and $J$ is invariant the relative density of $J$ in
  $\hat A$ is 1. The Julia set is closed, hence $\hat A\subset J$. This
  is impossible.
\end{plm}
 
\AuthorRefNames [MdMv2]

\end{document}